\documentclass[12pt]{amsart}
\usepackage{amssymb,amsmath}
\usepackage{xypic}

\theoremstyle{plain}
\newtheorem{lem}{Lemma}[section]
\newtheorem{pro}[lem]{Proposition}

\newtheorem{thm}[lem]{Theorem}
\newtheorem{cor}[lem]{Corollary}

\theoremstyle{definition}
\newtheorem{defn}[lem]{Definition}
\newtheorem{eg}[lem]{Example}
\newtheorem{rem}[lem]{Remark}

\newcommand{\bb}[1]{\mathbb{#1}}
\newcommand{\cl}[1]{\mathcal{#1}}
\newcommand{\mf}[1]{\mathfrak{#1}}

\newcommand{\eval}{\operatorname{eval}}
\newcommand{\Hol}{\operatorname{Hol}}
\newcommand{\Mul}{\operatorname{Mul}}
\newcommand{\range}{\operatorname{range}}
\newcommand{\spr}{\operatorname{spr}}

\renewcommand{\phi}{\varphi}
\newcommand{\ep}{\varepsilon}

\newcommand{\qand}{\quad\text{and}\quad}
\newcommand{\qfor}{\quad\text{for}\quad}

\newcommand{\ol}[1]{\overline{#1}}
\newcommand{\bsl}{\setminus}

\newenvironment{sbmatrix}{\left[\begin{smallmatrix}}%
{\end{smallmatrix}\right]}
\newenvironment{spmatrix}{\left(\begin{smallmatrix}}%
{\end{smallmatrix}\right)}

\begin{document}

\title{Generalized Berezin Transform and Commutator Ideals}
\author[K.R.Davidson]{Kenneth R. Davidson}
\thanks{First author partially supported by an NSERC grant}
\address{Pure Math.\ Dept.\\U. Waterloo\\Waterloo, ON\;
N2L--3G1\\CANADA}
\email{krdavids@uwaterloo.ca}
\author[R.G.Douglas]{Ronald G. Douglas}
\thanks{This research  was started while the second author was
visiting the Fields Institute on  research  leave from Texas A\&M
University and continued during  visits by him to the  Laboratoire de
Math\'ematique Pures, Universit\'e Bordeaux and the  Institut des Hautes
\'Etudes Scientifique.}
\address{Math.\ Dept.\\Texas A\&M U.\\
College Station, TX\ 77843-3368 }
\email{rgd@tamu.edu}
\date{}
\begin{abstract}
For a quasi-free module over a function algebra $A(\Omega)$, we define
an analogue of the Berezin transform and relate this to the quotient of
the C*-algebra it generates modulo the commutator ideal.
\end{abstract}
\maketitle

Certainly the best understood commutative Banach algebras are those 
that consist  of all the continuous complex-valued functions on a
compact Hausdorff space.  Indeed, most self-adjoint phenomena
involving them have been thoroughly  investigated. In particular, the
study of their representation theory as  operators on a Hilbert space,
which is essentially the spectral theory for  normal operators, shows
that such representations are defined by multiplication  on
$L^2$-spaces. Over the past few decades, other classes of operators
have been  introduced that are defined by functions, but which involve
more complicated  methods. One example is Toeplitz operators while
another example is the class of  pseudo-differential operators. In
both cases, one shows that the operators so  defined behave like the
functions used to define them, up to operators of lower  order.

To be more precise, let $H^2(\bb D)$ be the Hardy space of functions 
in $L^2(\bb  T)$  consisting of the functions with zero negative
Fourier coefficients and $P$ be  the projection of $L^2(\bb T)$ onto
$H^2(\bb D)$. The Toeplitz operator $T_\phi$ for the 
function $\phi$ in $L^\infty(\bb T)$ is defined on $H^2(\bb D)$ to be 
pointwise multiplication by $\phi$ followed by $P$. If $\mf T$ denotes the 
C*-algebra generated by the Toeplitz operators $T_\phi$ for $\phi$ a 
continuous function on $\bb T$, then $\mf T$ contains the algebra of 
compact  operators $\mf K$ and the quotient algebra $\mf T/\mf K$ is
isometrically  isomorphic to $C(\bb T)$. If instead of $H^2(\bb D)$ we
take the Bergman space  $B^2(\bb D)$ consisting of the functions in
$L^2(\bb D)$ which are a.e.\ equal  to a holomorphic function on $\bb
D$, then we can define ``Toeplitz operators''  analogously, and the
C*-algebra $\mf T'$ generated by the operators defined  by continuous
functions on the closed unit disk, contains $\mf K$ and the  quotient
algebra $\mf T'/\mf K$ is again isometrically isomorphic to $C(\bb T)$.

The description of these two examples can be carried over to the case 
of several  variables in more than one way. First, consider the
boundary $\partial \bb B^n$  of the unit ball $\bb B^n$ in $\bb C^n$
and the Hardy space of functions in 
$L^2(\partial \bb B^n)$ that have holomorphic extensions to the ball. 
Again, we  can define Toeplitz operators using functions defined on
$\partial\bb B^n$ and  the C*-algebra generated by the ones defined by
continuous functions again  contains the compact operators and the
quotient algebra is isometrically  isomorphic to $C(\partial \bb
B^n)$. Second, for the corresponding Bergman space  we can define
Toeplitz operators for functions on the closed ball and the 
C*-algebra generated by the continuous functions contains the compact 
operators and again, the quotient algebra is isometrically isomorphic
to $C(\partial\bb B^n)$. Thirdly, if we consider the polydisk $\bb D^n$
and the  Hardy subspace of $L^2(\bb T^n)$, then although the
C*-algebra generated by  the Toeplitz operators defined by continuous
functions on $\bb T^n$ contains the  ideal of   compact operators $\mf
K$, in this case $\mf K$ is not the commutator ideal $\mf  C$ in $\mf
T$. But $\mf C$ is proper and the quotient algebra is isometrically 
isomorphic to $C(\bb T^n)$.

The phenomenon we want to consider in this note concerns an algebra of 
operators  defined by functions, the C*-algebra that the algebra
generates and the  quotient algebra defined modulo the commutator
ideal. The results we obtain will  contain all the above examples and
show that, in a certain sense, it is the  domain 
that is important, rather than the particular Hilbert space on which 
the  operators are defined. We will accomplish this by generalizing
the notion of a  transform introduced by Berezin \cite{B} in
connection with quantization. In our  approach, the Hilbert space is
closely related to the algebra of holomorphic  functions on the domain
and involves a notion of kernel Hilbert space that is a  module over
the algebra.

We build on results of several authors who have studied the Berezin 
transform.  In particular, there is the result of McDonald and
Sundberg \cite{M-S}  concerning the 
nature of the transform for the C*-algebra of Toeplitz operators 
defined by  bounded holomorphic multipliers on the Bergman space.  
Some of our  proofs are closely related to arguments in \cite{M-S}.
Further, some issues we  study overlap questions raised by Arveson in
\cite{Arv}. We also relate a  further 
study of Sundberg \cite{S} on the C*-algebras generated by 
Toeplitz-like  operators to our results including a relation between
the property of a module  multiplier having closed range to the
behavior of its Gelfand transform.  Finally, we study the 
relationship between the kernel ideal of  the generalized Berezin
transform and the commutator ideal.

\section{Definition and Basic Properties}\label{sec1}

Regardless of the original motivation of Berezin for introducing it 
(cf.\ \cite{B}), the Berezin transform essentially provides a kind of
``symbol'' for  certain natural operators on Hilbert spaces of
holomorphic functions. For our  generalized transform, we will use
kernel Hilbert spaces over bounded domains in 
$\bb C^n$, which are also contractive Hilbert  modules for natural 
function  algebras over the domain. More precisely, we use the concept
introduced in  \cite{D-M} for the study of module resolutions.

For $\Omega$ a bounded domain in $\bb C^n$, let $A(\Omega)$ be the 
function  algebra obtained as the completion of the set of functions
that are holomorphic  in some neighborhood of the closure of $\Omega$.
For $\Omega$ the unit ball $\bb  B^n$ or the polydisk $\bb D^n$ in
$\bb C^n$, we obtain the familiar ball and  polydisk algebras $A(\bb
B^n)$ and $A(\bb D^n)$, respectively. The Hilbert space 
$\cl M$ is said to be a {\em contractive Hilbert module\/} over 
$A(\Omega)$ if $\cl M$ is a unital  module over $A(\Omega)$ with 
module map $A(\Omega)\times \cl M\to \cl M$ such that
\[
 \|\phi f\|_{\cl M} \le \|\phi\|_{A(\Omega)} \|f\|_{\cl M}
 \qfor \phi \in A(\Omega) \qand f \in \cl M.
\]
The space $\cl R$ is said to be a {\em quasi-free Hilbert module of 
rank\/} $m$, $1\le m\le \infty$, over $A(\Omega)$, if it is obtained
as the completion of an inner  product on the algebraic tensor product
$A(\Omega) \otimes \ell^2_m$ such that
\begin{enumerate}

\item $\eval_{\pmb{z}}: A(\Omega) \otimes \ell^2_m \to \ell^2_m$,
 the canonical evaluation map at the point $\pmb{z}$, is bounded
 for $\pmb{z}$ in $\Omega$ and locally uniformly bounded on
 $\Omega$;

\item $\|\Sigma\phi\theta_i \otimes x_i\|_{\cl R} \le 
 \|\phi\|_{A(\Omega)} \, \| \Sigma \theta_i \otimes x_i\|_{\cl R}$
 for $\phi$, $\{\theta_i\}$ in $A(\Omega)$ and $\{x_i\}$ in
 $\ell^2_m$; and

\item for $\{F_i\}$ a sequence in $A(\Omega)\otimes \ell^2_m$ Cauchy
 in the $\cl R$-norm, it follows that $\eval_{\pmb{z}}(F_i)\to 0$ for
 all $\pmb{z}$  in $\Omega$ if and only if $\|F_i\|_{\cl R}\to 0$.

\end{enumerate}
Here, $\ell^2_m$ is the $m$-dimensional Hilbert space.

In \cite{D-M}, another characterization and other properties of 
quasi-free  Hilbert modules are given. This concept is closely related
to the notions of sharp and  generalized Bergman kernels studied by
Curto and Salinas \cite{C-S}, Agrawal and Salinas \cite{A-S} and
Salinas \cite{salinas}.

Although much of what follows is valid for the case of $m=\infty$, we
confine  our attention here to the finite rank case. Hence, the
evaluation function, $\eval_{\pmb{z}}: \cl R\to \ell^2_m$, which is
defined to be the  extension of evaluation on
$A(\Omega)\otimes\ell^2_m$, is onto and an $\cl  R$-valued holomorphic
function in $\pmb{z}$.  If $\eval^*_{\pmb{z}}: \ell^2_m\to \cl R$ is
the operator adjoint, then it is bounded below and is 
anti-holomorphic in $\pmb{z}$. Moreover, the lower bound is continuous
in $\pmb{z}$ since $\eval_{\pmb{z}} \eval^*_{\pmb{z}}$ is a
real-analytic  function in $\pmb{z}$.

The following lemma provides some well-known connections between 
$\eval_{\pmb{z}}$, module multiplication, and the kernel function. 
Let $M_\phi$ denote the operator on $\cl R$ defined by  module
multiplication.

\begin{lem}\label{eval}
For $\pmb{z}$ in $\Omega$ and $\phi$ in $A(\Omega)$, one has 
\[
 M^*_\phi \eval^*_{\pmb{z}} =\overline\phi(\pmb{z}) \eval^*_{\pmb{z}} 
 \qand 
 \eval_{\pmb{z}} M_\phi = \phi(\pmb{z}) \eval_{\pmb{z}} .
\]
Moreover, 
\[
 K(\pmb{w},\pmb{z}) = \eval_{\pmb{w}} \eval^*_{\pmb{z}}:
 \Omega\times \Omega\to \mf L (\ell^2_m)
\]
is the kernel function for $\cl R$.
\end{lem}

\begin{proof}
The two identities in the first sentence follow from the fact that 
\[
 (\eval_{\pmb{z}} M_\phi)f = 
 \eval_{\pmb{z}}(\phi f) = 
 \phi(\pmb{z}) f(\pmb{z}) = 
 \phi(\pmb{z})  \eval_{\pmb{z}} f
\]
for $\phi \in A(\Omega)$, $f \in \cl R$ and $\pmb{z} \in \Omega$.

Finally, for $\pmb{w}$ and $\pmb{z}$ in $\Omega$ and $x$ in $\ell^2_m$, the 
function $f(\pmb{w}) = 
 K(\pmb{w},\pmb{z})x$ is in $\cl R$ and 
\[
 \langle g,f\rangle_{\cl R} =
 \langle g, K(\cdot,\pmb{z})x\rangle_{\cl R}=
 \langle g(\pmb{z}),x\rangle_{\ell^2_m}  
 \qfor g \in \cl R
\]
by the definition of the kernel function $K$. 
Since $g(\pmb{z}) = \eval_{\pmb{z}}g$, we 
have 
\[
 \langle g(\pmb{z}),x\rangle =
 \langle\eval_{\pmb{z}}g,x\rangle =
 \langle g, \eval^*_{\pmb{z}} x\rangle_{\cl R}
\]
and therefore, 
$K(\cdot,\pmb{z})x = \eval^*_{\pmb{z}}x$. 
Applying $\eval_{\pmb{w}}$ to both sides yields $K(\pmb{w},\pmb{z}) x=
\eval_{\pmb{w}} \eval^*_{\pmb{z}}x$ for $x$ in $\ell^2_m$ which is the
desired result.
\end{proof}

To define the transform we need the polar form of evaluation on $\cl R$. Thus, 
set $\eval^*_{\pmb{z}} = V_{\pmb{z}}Q_{\pmb{z}}$, where 
1)~~$V_{\pmb{z}}$ is an isometry from 
$\ell^2_m$ into $\cl R$ which is real-analytic in $\pmb{z}$ and 
2)~~$Q_{\pmb{z}}$ is a positive, invertible operator on $\mf
L(\ell^2_m)$, which is also real-analytic  in $\pmb{z}$. The
properties of $\eval_{\pmb{z}}$ in Lemma~\ref{eval} carry over 
to $V_{\pmb{z}}$ to yield the following results:

\begin{lem}\label{eval2}
For $\phi$ in $A(\Omega)$ and $\pmb{z}$ in $\Omega$, it follows that
\begin{enumerate}
 \item $M^*_\phi V_{\pmb{z}} = \overline{\phi(\pmb{z})}V_{\pmb{z}}$ and
 \item $V^*_{\pmb{z}}M_\phi = \phi(\pmb{z}) V^*_{\pmb{z}}$.
\end{enumerate}
\end{lem}

We can now define the Generalized Berezin Transform (GBT) on
$\cl R$.

\begin{defn}\label{def_GBT}
If $\cl R$ is a quasi-free Hilbert module of finite 
rank over $A(\Omega)$,
then the {\em Generalized Berezin Transform\/} $\Gamma$ maps   
$\mf L(\cl R)$ to $\mathrm{C}_b(\Omega, \mf L(\ell^2_m))$,
the space of bounded continuous $m \times m$ matrix-valued
functions on
$\Omega$, by $(\Gamma X)(\pmb{z}) =  V^*_{\pmb{z}}XV_{\pmb{z}}$
for $X$ in $\mf L(\cl R)$ and $\pmb{z}$ in $\Omega$.
\end{defn}

\begin{lem}\label{lem_GBT}
The GBT has the following elementary properties:
\begin{enumerate}
 
\item $\Gamma$ is contractive:
 $\|(\Gamma X)(\pmb{z})\|_{\mf L(\ell^2_m)} \le \|X\|_{\mf L(\cl R)}$ 
 for $X \in \mf L(\cl R)$ and $\pmb{z} \in \Omega$,

\item $\Gamma$ is linear: 
 $\Gamma(a_1X_1 + a_2X_2) = a_1\Gamma(X_1) + a_2\Gamma(X_2)$ for 
 $a_1,a_2 \in \bb C$, $X_1,X_2 \in \mf L(\cl R)$, and

\item $\Gamma$ is self-adjoint:
 $(\Gamma X)(\pmb{z})^*  = \Gamma(X^*)(\pmb{z})$
 for $X \in \mf L(\cl R)$ and $\pmb{z} \in \Omega$.

\item $\Gamma(X)(\pmb z)$ is continuous on $\Omega$.

\end{enumerate}
\end{lem}

\begin{proof} Only (4) requires any comment.
The map sending $\pmb{z}$ to $\eval_{\pmb{z}}$ is continuous on
$\Omega$, and $\eval_{\pmb{z}}$ is surjective for all $\pmb{z} \in
\Omega$.  It follows that  the functions taking $\pmb{z}$ to
$V_{\pmb{z}}$ and $V^*_{\pmb{z}}$ are also continuous. 
Therefore, the map $\Gamma(X)(\pmb{z}) = V^*_{\pmb{z}} 
XV_{\pmb{z}}$ is a bounded, continuous  
$\mf L(\ell^2_m)$-valued function on $\Omega$. 
\end{proof}

If $V_{\pmb{z}}$ is any isometry-valued function from $\Omega$ to 
$\mf L(\ell^2_m,\cl R)$, one could define a GBT having the properties
in this lemma but little more. The  interest in our setting comes
about because of the relationship between the GBT  and the operators
defined by the algebra $A(\Omega)$. In general, one is  interested in
the behavior of the GBT on the C*-algebra generated by the  operators
defined by module multiplication or the corresponding algebra
obtained  when $A(\Omega)$ is replaced by a larger, natural function
algebra of bounded  holomorphic functions (cf.\ Section~\ref{sec4}).

The following proposition is easily deduced from
Lemma~\ref{eval2}.

\begin{pro}\label{prop_GBT}
For $\phi,\psi$ in $A(\Omega)$ and $X$ in $\mf L(\cl R)$ it follows 
that:
\begin{enumerate}

\item $\Gamma(M_\phi X)(\pmb{z}) = \phi(\pmb{z})(\Gamma X)(\pmb{z})$ 
 and $\Gamma(XM^*_\phi)(\pmb{z}) =
 \overline\phi(\pmb{z})(\Gamma X)(\pmb{z})$ 
 for $\pmb{z}$ in $\Omega$,

\item $\Gamma(M_\phi) = \phi$ and $\Gamma(M^*_\phi) = 
 \overline{\phi}$, and

\item $\Gamma(M_\phi M^*_\psi) = \phi \overline{\psi}$.

\end{enumerate}
\end{pro}

These results show that the ``symbol'' defined by the GBT agrees with 
the  multiplier for operators defined by module multiplication for
functions in $A(\Omega)$ or their complex conjugates. Also, it is
multiplicative on products  of the form $\phi\overline\psi$ for
$\phi,\psi$ in $A(\Omega)$; but it is  easy to check that
$\Gamma(M_\phi M^*_\psi)\ne \Gamma(M^*_\psi  M_\phi)$, in general.
This comes down to evaluating $\Gamma$ on commutators  of the form
$[M_\phi, M^*_\psi]$ for $\phi,\psi$ in $A(\Omega)$.

If one considers the classical examples in which $\cl R$ is either the
Hardy or Bergman  modules on the unit disk $\bb D$, one  sees that the
Berezin transform is  nicely-behaved for the C*-algebra generated by
module multipliers {\em  only\/} on the boundary of 
$\bb D$ and vanishes there on the commutator ideal. More interesting is the 
behavior of the Berezin Transform on the corresponding algebra defined by  
module multiplication by the algebra $H^\infty(\bb D)$ of bounded
holomorphic functions  on $\bb D$ which we'll consider in
Section~\ref{sec4}.

To get started we consider the basic relationship between the 
operators defined  by module multiplication and the quotient algebra
defined by the commutator  ideal. Understanding this quotient algebra
is equivalent to characterizing the  multiplicative linear 
functionals on the algebra.

Let $\mf T(\cl R)$ denote the C*-algebra generated by 
$\{M_\phi : \phi\in A(\Omega)\}$ acting on $\cl R$ and $\mf C(\cl R)$ 
be  the closed two-sided ideal in $\mf T(\cl R)$ generated by the
commutators in $\mf T(\cl R)$. Then $\mf T(\cl R)/\mf C(\cl R)$ is a
commutative C*-algebra,  possibly (0), and hence is isometrically
isomorphic to $C(X_{\cl R})$ for some compact Hausdorff space  
$X_{\cl R}$.  Thus the points of $X_{\cl R}$ are precisely the
characters on $\mf T(\cl R)$.  So we call $X_{\cl R}$ the
\textit{character space of $\cl R$}. Our first result relates the
character space to the  maximal ideal space $M_A$ of $A(\Omega)$.

\begin{thm}\label{char_space}
For $\cl R$ a quasi-free Hilbert module of finite rank, the map taking
$\phi$ to  $M_\phi + \mf C(\cl R)$ extends to a surjective 
$*$-homomorphism $\tau$ from $C(M_A)$ to 
$\mf T(\cl R)/\mf C(\cl R) = C(X_{\cl R})$. 
Consequently, one can identify $X_{\cl R}$ with a closed subset of
$M_A$.
\end{thm}

\begin{proof}
If $\mf C(\cl R) = \mf T(\cl R)$ whence $X_{\cl R}$ is empty, then 
$\tau$ is the  zero  map and the result is vacuous; hence we assume
that $X_{\cl R}$ is non-empty.  Also, we abbreviate $A = A(\Omega)$ in
the following proof.

Consider the diagram
\[
\xymatrix{
 0 \ar[r] & \mf C(\cl R) \ar[r]^i & \mf T(\cl R) \ar[r]^q
 & C(X_{\cl R}) \ar[r] & 0
 \\ 
 && A(\Omega) \ar[u]^\mu \ar[r] \ar[ur]^{\hat{\tau}}
 & C(M_A) \ar[u]_\tau }
\]
where the map $\mu$ from $A(\Omega)$ to $\mf T(A)$ sends $\phi$ to
$M_\phi$, the map from $A(\Omega)$ to  $C(M_A)$ is the Gelfand map, 
and $\hat{\tau} = q \mu$ is defined by composition. Then 
$\tau$ exists by the universal properties of the Gelfand map. 
Moreover, since  the image of $A(\Omega)$ in $\mf T(A)$ generates it,
we have that the range of $\hat{\tau}$ generates $C(X_{\cl R})$.
Therefore, the range of $\tau$  generates $C(X_{\cl R})$ and hence
separates the points of $X_{\cl R}$ which  implies that we can
identify $X_{\cl R}$ as a closed subset of $M_A$.
\end{proof}

As we mentioned above, for  the Hardy and Bergman modules, the 
quotient  algebras $\mf T(H^2(\bb D))/ \mf C(H^2(\bb D))$ and    
$\mf T(B^2(\bb D))/\mf  C(B^2(\bb D))$ are both  equal to
$C(\partial\bb D)$. Hence,  for both $\cl R = H^2(\bb D)$ and
$B^2(\bb D)$, the character space is equal to the \v Silov boundary.
We will see in Example~\ref{wtd_shift} that this is not always the
case. An even more dramatic example will be given in
Example~\ref{simple} in which we show for a module which is similar
to the Hardy module (but not contractive) that the character space
can be empty. This is in striking contrast to the fact that for
contractive modules over nice domains including the unit disk, the
character space must contain the whole \v Silov boundary. 
See Corollary~\ref{Shilov}.

\section{Pointed Function Algebras}\label{sec2}

We introduce a property for points in the maximal ideal space of a
function algebra which is related to the well-known notions of
peak point and p-point.  We shall give a few cases where this
property can be explicitly exhibited.  Once this work was close to
completion, we asked John Wermer if he could establish this
property in certain examples.  The result is a complete answer by
Izzo and Wermer cited below that reduces our notion to the
classical ones.

\begin{defn}\label{defn_pointed}
The function algebra $A$ is said to be {\em pointed} at $\alpha$ in
$M_A$ if $A$ is the closed linear span of the set 
$\{ \phi \in A : |\phi(\alpha)| = \|\phi\|_A \}.$

A function algebra $A$ is said to be {\em pointed\/} if it is pointed
on a dense subset of the \v Silov boundary $\partial A$ of $A$.

We will say that a domain $\Omega$ is pointed if $A(\Omega)$ is
pointed.
\end{defn}

Note that since the set of functions achieving their maxima at a
point $\alpha$ is closed under  multiplication, the  closed linear
span is automatically a closed subalgebra of $A$.

Recall that a point $\alpha$ in $M_A$ is a \textit{p-point} if for 
every neighbourhood $U$ of $\alpha$, there is a function $\phi$ in
$A$ which attains its maximum modulus at $\alpha$ and nowhere
outside of $U$.  It is a \textit{peak point} if there is a function
$\phi$ which attains its maximum precisely at $\alpha$.  In the
metric case, it is well-known that p-points are peak points.  

The following easy lemma shows that being pointed implies being a
p-point.

\begin{lem}\label{peak}
If function algebra $A$ is pointed at $\alpha$, then $\alpha$ is a
p-point.  In the metrizable case, it is a peak point.
\end{lem}

\begin{proof}
If $U$ is a neighborhood of $\alpha$, for each point $\beta$ in $M_A
\setminus U$, choose a function $f$ in $A$ such that 
$\| f \| = f(\alpha) = 1$ and $f(\beta) \ne 1$. It is easy to modify
$f$ so that its range meets the unit circle only at 1. A finite
cover yields a finite family of such functions $f_i$, $1 \le i \le
k$ so that for each $\beta$ in $M_A \setminus U$, there is an $i$
with
$|f_i(\beta)| < 1 = f_i(\alpha) = \| f_i \|$.  The average of the
$f_i$'s peaks at $\alpha$ and has modulus less than 1 on $M_A
\setminus U$.  So $\alpha$ is a p-point.
\end{proof}

The converse was obtained by Izzo and Wermer \cite{I-W}.  It is a
short but clever argument, and yields a much stronger result than we
anticipated.

\begin{thm}[Izzo--Wermer]\label{IW}
If $\alpha$ is a p-point of a function algebra $A$, then $A$ is
pointed at $A$.  If $\alpha$ is a peak point, then every function in
$A$ is the linear combination of two functions which peak precisely
at $\alpha$.
\end{thm}

In the metric case, the peak points form a dense $G_\delta$ of the \v
Silov boundary. In general, the p-points are dense (but perhaps not
even Borel).  See Gamelin \cite{Gam}. Thus we obtain:

\begin{cor}\label{pointed}
Every function algebra is pointed.
\end{cor}

Even though this general result eliminates the need to demonstrate
that particular function algebras are pointed, it is of interest to
see when this can be achieved in practice.  So we provide a few
examples.

\begin{eg}\label{Jordan}
Consider finitely connected domains in $\bb C$ whose boundary 
consists  of a finite union of Jordan curves. 
Let us begin with the case of the annulus 
$\bb A  = \{z \in \bb C : 0<r<|z| <  R < \infty\}$. 
Fix $\alpha$ on the outer circle. 
The functions $\{z^n\}$, $n\ge 0$, satisfy the 
requirement that $\|z^n\| =  R^n = |\alpha^n|$. Moreover, the function 
$\phi_\lambda(z) = z + \lambda z^{-n}$, for $n>0$ and $\lambda$ in 
$\bb  C$, takes its maximum on either the inner or outer circles. For
$\lambda$ of  sufficiently small absolute value but non-zero, the
maximum is achieved on the  outer circle. By choosing the phase of
$\lambda$ correctly, call it $\lambda_n$,  we can force the maximum
absolute value of $\phi_{\lambda_n}$ to occur at $\alpha$. Since the
closed linear span of the set of functions 
$\{z^n, z + \lambda_n z^{-n} \}$ is $A(\bb A)$, 
we see that $\bb A$ is pointed at each $\alpha$ on the outer circle. 
But being  pointed is a conformal invariant in the sense that if
$\Omega$ is pointed at  $\alpha$ in $\partial\Omega$ and $\mu$ is a
conformal self-map on $\Omega$ that  extends to $\partial\Omega$,
then $\Omega$ is also pointed at $\mu(\alpha)$.  Therefore, $\bb A$
is a pointed domain as are all one-connected domains with  Jordan
curves as boundaries.

Now consider an arbitrary bounded finitely connected domain $\Omega$ 
in $\bb C$  with boundary consisting of Jordan curves. By making a
conformal transformation  on $\Omega$ which extends to
$\partial\Omega$, we can assume the outer boundary  is a circle. Let
$\alpha$ be a point on the bounding outer circle of $\Omega$.  If 
$z_0$ is a point inside this circle not in $\Omega$, we can repeat the 
above argument to exhibit a function $z + 
\lambda_n (z-z_0)^{-n}$ which peaks at $\alpha$. Doing this for a 
point in each bounded component of $\bb C\setminus \Omega$, we obtain 
a set of  functions in $A(\Omega)$ each of which peaks at $\alpha$ and
for which their  closed linear span equals $A(\Omega)$. Therefore,
$\Omega$ is a pointed domain.
\end{eg}

\begin{eg}\label{product}
An easy argument shows that the product
$\Omega_1\times \Omega_2$ in $\bb C^{m+n}$ of two pointed domains
$\Omega_1$ in $\bb C^m$ and $\Omega_2$ in  $\bb C^n$ is also pointed.
The \v Silov boundary of the product is the product of the \v Silov boundaries.
The key observations are 
(1)~~that the closed linear  span of
functions of the form $\phi\psi$ for $\phi$ in $A(\Omega_1)$ and 
$A(\Omega_2)$ generates $A(\Omega_1\times \Omega_2)$ and 
(2)~~that $|\phi(\alpha_1)\psi(\alpha_2)| =
\|\phi\psi\|_{A(\Omega_1\times \Omega_2)}$ if $|\phi(\alpha_1)| =
\|\phi\|_{A(\Omega_1)}$ and 
$|\psi(\alpha_2)| = \|\psi\|_{A(\Omega_2)}$. 
Therefore, the polydisk is pointed.
\end{eg}

The following argument is not nearly as definitive as the Izzo--Wermer
result, but provides a rather different approach.

\begin{pro}\label{ball}
The ball algebra $A(\bb B^n)$ is pointed.
\end{pro}

\begin{proof} Let $n=2$.
Observe that the following functions peak at $(1,0)$: $1$, $z_1$,
$z^2_1 + z^2_2$ and $z^2_1+z^3_2$. In addition, let $h\in A(\bb D)$
with $\|h\|_\infty \le 1$ and consider 
$f(z_1,z_2) = (1+z_1) z_1/2 + (1-z_1) z_2 h(z_1)/2$.  Then by
Cauchy--Schwartz,
\begin{align*}
 |f(z_1,z_2)|^2 &\le \left(\left| \frac{1+z_1}2\right|^2 + 
 \left|\frac{1-z_1}2\right|^2\right) (|z_1|^2 +  |z_2|^2 \|h\|^2)\\
 &\le \frac{1+|z_1|^2}2 \le 1.
\end{align*}
So this function also peaks at (1,0).

The span of all functions peaking at (1,0) is an algebra. So this span 
includes  all polynomials of the form $\sum_i h_i(z_1) z^i_2$ provided
that $h_1(1) = 0$.  But the functional taking a polynomial to
$h_1(1,0)$, namely $\phi(f) = \frac\partial{\partial z_2} f(1,0)$ is
discontinuous in the ball norm. Thus the  closure in $A(\bb B^2)$ is
the whole algebra.

This argument readily generalizes to $n>2$.
\end{proof}

\begin{eg}\label{convex}
It is easy to see that if $\Omega$ is contained in a large ball which 
is tangent to $\Omega$ at a point $\alpha$, then $\Omega$ is pointed
there. So any convex domain which has strictly positive curvature at
each point on the boundary will be a pointed domain.
\end{eg}

\section{The Commutator Ideal}\label{sec3}

We now consider the nature of the range of the Generalized Berezin
Transform.  Our goal is to establish that the quotient of $\mf T(\cl
R)$ by its commutator ideal is large, and our tool for establishing
that $\mf T(\cl R)$ has many characters is the GBT.

By Proposition~\ref{prop_GBT}, the map taking $\pmb{z}$ to
$\Gamma X (\pmb{z}) = V^*_{\pmb{z}}  XV_{\pmb{z}}$ is a
bounded, continuous $\mf L(\ell^2_m)$-valued function on $\Omega$
for $X$ in $\mf L(\cl R)$. Hence we can extend $\Gamma X$ to the 
Stone--\v Cech  compactification $\beta\Omega$ of $\Omega$.
The inclusion of $A(\Omega)$ in $C(\beta\Omega)$ defines a map 
$\rho: \beta\Omega\to M_A$. Suppose for convenience that $m=1$.
For $\alpha$ in $\beta\Omega$, we can ask
if the state on  $\mf T(\cl R)$ defined by $X\to \Gamma(X)(\alpha)$
is multiplicative. If it is,  then it must coincide with the extension
of the map on $\mf T(\cl R)$ defined by 
$M_\phi\to \widehat\phi(\rho(\alpha))$. 
Thus there is a close relationship between the Berezin transform on
$\mf T(\cl R)$ and the Gelfand transform on 
$\mf T(\cl R)/\mf C(\cl R)$, which we investigate in this section.
Our approach is based on the fact that $A(\Omega)$ is pointed. 

The following lemma is the key observation that takes advantage
of functions which attain their maximum modulus at $\alpha$.
Let $P_{\pmb{z}} = V_{\pmb{z}} V_{\pmb{z}}^*$ be the projection onto
the range of R$V_{\pmb{z}}$ for ${\pmb{z}} \in \Omega$.

\begin{lem}\label{asymp_reducing}
Suppose that $\alpha \in M_A$ and $\phi$ is a function in $A(\Omega)$
such that  $|\widehat\phi(\alpha)| = \|\phi\|_{A(\Omega)}$.
Then $\lim\limits_{\pmb{z}\to \alpha} \| [P_{\pmb{z}}, M_\phi] \| =
0$.  That is, the range of $V_{\pmb{z}}$ is asymptotically reducing
for $M_\phi$ as ${\pmb{z}}$ approaches $\alpha$.

For any $\psi \in A(\Omega)$ and $S \in \mf L(\cl R)$,  it follows
that 
\[
 \lim\limits_{\pmb{z}\to \alpha}
 \|\Gamma(S [M^*_\psi,M_\phi]) (\pmb{z})\|_{\mf L(\ell^2_m)} = 0 .
\]
\end{lem}

\begin{proof}
Decompose $\cl R = (\range V_{\pmb{z}})\oplus (\range  V_{\pmb{z}})^\perp$,
and write $M_\phi$ as a matrix relative to the decomposition. 
Since $ V_{\pmb{z}}^* M_\phi= \phi(\pmb{z}) V_{\pmb{z}}^*$,
we obtain
\[
 M_\phi\simeq
\begin{bmatrix}
 \phi(\pmb{z}) I_{\range V_z} & 0 \\
 C_{\pmb{z}} & D_{\pmb{z}} 
\end{bmatrix} 
\qand
 P_{\pmb{z}} = 
\begin{bmatrix}
 I_{\range V_z} &0 \\
 0 & 0
\end{bmatrix}.
\]
Moreover, since
\[
 \|M_\phi\|_{\mf L(\cl R)} \le \|\phi\|_{A(\Omega)} = 1 \qand  
 |\phi(z)| \le \|\phi\|_{A(\Omega)} = 1,
\]
it follows that $C_{\pmb{z}} = (1-|\phi(z)|^2)^{1/2} C'_{\pmb{z}}$, 
where $\|C'_{\pmb{z}}\|\le 1$, whence we obtain that
$\| C_{\pmb{z}}\| \le (1-|\phi(z)|^2)^{1/2}$.
Therefore
\[
 \| [M_\phi, P_{\pmb{z}}] \| =
 \Big\| 
 \begin{bmatrix} 0 & 0 \\ C_{\pmb{z}} & 0 \end{bmatrix} 
 \Big\|
 \le (1-|\phi(z)|^2)^{1/2} .
\]

Let us write
\[
 M^*_\psi\simeq
\begin{bmatrix}
 \overline{\psi(\pmb{z})} I_{\range V_z} & X_{\pmb{z}} \\
 0 & Y_{\pmb{z}} 
\end{bmatrix} 
\qand
S = 
\begin{bmatrix}
 S_{11} & S_{12} \\
 S_{21} & S_{22}
\end{bmatrix} .
\]
A straightforward calculation now yields the $1,1$ entry of the
operator $S [M^*_\psi,M_\phi] $:
\[
 \Gamma(S [M^*_\psi,M_\phi] )(\pmb{z}) = V_{\pmb{z}}^*
 \big( S_{11} X_{\pmb{z}} + S_{12} (Y_{\pmb{z}} - \ol{\psi}(\pmb{z})
 \big) C_{\pmb{z}} V_{\pmb{z}}.
\]
This is easily bounded by 
$3 \| \psi \| \, \| S \| \, \| C_{\pmb{z}} \|$ which tends to 0 as
$\pmb{z}$ converges to $\alpha$.
\end{proof}

\begin{thm}\label{extend_peak}
If $\Omega$ is pointed at $\alpha \in M_A$, then
$
 \lim\limits_{\pmb{z}\to \alpha} \Gamma(T) (\pmb{z})
$
exists for all $T$ in $\mf T(\cl R)$, and this map defines a
character of $\mf T(\cl R)$ corresponding to evaluation at $\alpha$.
In particular,  the GBT evaluated at $\alpha$ in this way annihilates the
commutator ideal of $\mf T(\cl R)$.  
\end{thm}

\begin{proof}
By Lemma~\ref{asymp_reducing}, the ranges of $V_{\pmb{z}}$ are
asymptotically reducing for operators $M_\phi$ corresponding to those
functions $\phi$ in $A(\Omega)$ which peak at $\alpha$.
By hypothesis, these functions generate $A(\Omega)$.
However it is a straightforward calculation to see that the set of
operators $T \in \mf L( \cl R)$ with the property that
$
 \lim\limits_{\pmb{z}\to \alpha} \| [ T, P_{\pmb{z}} ] \| = 0
$
is a C*-algebra.  Consequently it contains all of $\mf T(\cl R)$.

The limit $\lim\limits_{\pmb{z}\to \alpha} \Gamma(M_\phi)(\pmb{z}) = 
\phi(\alpha) I_{\mf L(\ell^2_m)}$ exists for all $\phi \in
A(\Omega)$. Moreover, it follows that compression to the ranges of
$V_{\pmb{z}}$  is asymptotically multiplicative.
Therefore, this limit exists for all $T \in \mf T(\cl R)$, and
the map is a $*$-homomorphism.
As the generators are sent to a family of scalar operators, the
image consists of scalars.
Hence this determines a (unique) character extending evaluation at
the point $\alpha$. 
In particular, this map annihilates the commutator ideal.
\end{proof}

Arguments such as those given above have been used in studying the 
Berezin  transform on the Hardy and Bergman modules for $\bb D$ (cf.\
\cite{H-K-Z,M-S,A-Z}) but the role of pointedness was only implicit.  
A somewhat related notion appears in \cite{Arv} where one assumes that
the collection of unitary operators in the closure of the  set 
$\{M_\phi + M^*_\psi : \phi,\psi \in A(\Omega)\}$ generates $\mf
T(\cl R)$.
Our results hold without special hypotheses.

\begin{thm}\label{main_thm}
Let $\Omega$ be a domain in $\bb C^n$ and let $\cl R$ 
be a finite rank  quasi-free Hilbert module over $A(\Omega)$.
Then $\overline\Gamma = \Gamma|_{\mf  T(\cl R)}$ extends to a 
$*$-homomorphism from $\mf T(\cl R)$ to $C(\partial A)$.
The null space of $\overline\Gamma$, denoted $\gamma(\mf T(\cl R))$
and  called the Berezin nullity, is an ideal which contains
$\mf C(\cl R)$.
\end{thm}

\begin{proof}
The Izzo--Wermer Theorem~\ref{IW} shows that $A(\Omega)$ is pointed.
Hence by Theorem~\ref{extend_peak},
\[
 \overline\Gamma(T)(\alpha) :=
 \lim\limits_{\pmb{z}\to \alpha} \Gamma(T) (\pmb{z})
\]
is defined at all p-points in the \v Silov boundary.
That is, evaluation at $\alpha$ extends to a character on the
C*-algebra $\mf T(\cl R)$.

The set of characters on $\mf T(\cl R)$ is a compact set $X_{\cl R}$.
Restriction to $A(\Omega)$ yields a map of $X_{\cl R}$ into $M_A$.
This map is injective because $\{ M_\phi : \phi \in A(\Omega) \}$
generated $\mf T(\cl R)$.  
From the previous paragraph, $X_{\cl R}$ contains all p-points, and
hence contains their closure, which is the \v Silov boundary.

It follows that $\overline\Gamma(T)$ extends by continuity on $\partial A$
to all points of the \v Silov boundary, yielding a $*$-homorphism of 
$\mf T(\cl R)$ into $C(\partial A)$.
Since $A(\Omega)$ separates points in $M_A$ and thus in $\partial A$,
this map is surjective by the Stone--Weierstrass theorem.
\end{proof}

In many situations which occur in practice, every point in the
topological boundary of $\Omega$ is a peak point.  In this case, a
much stronger conclusion is possible.

\begin{cor}\label{all_points}
Assume that the \v Silov boundary $\partial A$ of $A(\Omega)$ coincides with
the topological boundary $\partial \Omega$, and furthermore that
every point in $\partial A$ is a peak point.
Then the generalized Berezin transform $\Gamma(T)$ is uniformly
continuous on
$\Omega$ and extends to $\overline{\Omega}$ and takes scalar
values on the boundary.

If every point of $\partial A$ is a peak point but $\partial A$ does not coincide
with the topological boundary $\partial \Omega$, then the GBT
extends to be continuous on $\Omega \cup \partial A$.
\end{cor}

\begin{proof} It is routine to verify that once we know (from Theorem~\ref{main_thm})
that $\Gamma(T)(\pmb{z})$ extends by continuity to every point on
$\partial\Omega$, that $\Gamma(T)$ is continuous on
$\overline{\Omega}$.  Moreover since evaluation at each boundary
point is a character, it takes scalar values there.

Likewise, if we only know continuity of the GBT at each popint of $\partial A$,
then we at least obtain continuity on $\Omega \cup \partial A$.  But as
we then have no information about the GBT as it approaches boundary points which
are not in the \v Silov boundary, we cannot make any claim to uniform continuity.
\end{proof}

\begin{rem}
We have not been able to determine whether or not the GBT extends by continuity to
the \v Silov boundary when there are non-peak points.  The well-known examples
all fall within the purview of the previous corollary, and so do
not provide any insight into this question.
\end{rem}

What we do obtain in general is the fact that the character space always contains
the \v Silov boundary.  Consider the diagram:
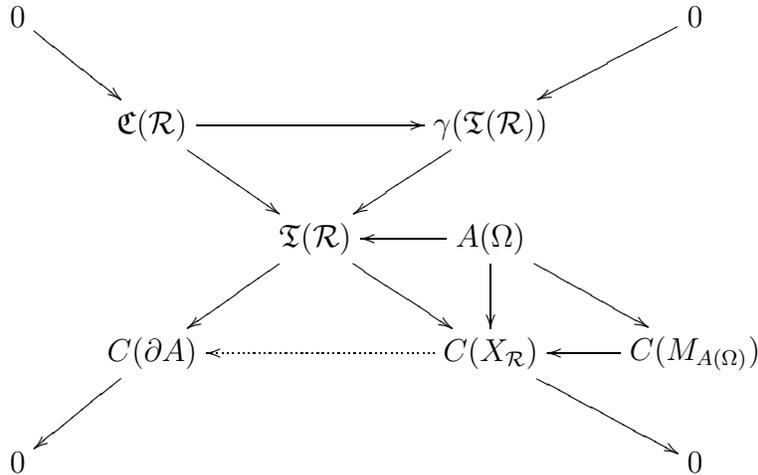
\[
\xymatrix{
 0 \ar[dr] &&&& 0 \ar[dl] \\
 & \mf C(\cl R) \ar[rr] \ar[dr] && \gamma(\mf T(\cl R)) \ar[dl] \\
 && \mf T(\cl R) \ar[dl] \ar[dr] & A(\Omega) \ar[l] \ar[d] \ar[dr] \\
 & C(\partial A) \ar[dl] && C(X_{\cl R}) \ar[dr] \ar@{.>}[ll] &
 C(M_{A(\Omega)}) \ar[l] \\
 0 &&&& 0
}
\]

\begin{cor}\label{Shilov}
Let $\Omega$ be a domain in $\bb C^n$ and let $\cl R$ 
be a finite rank  quasi-free Hilbert module over $A(\Omega)$.
Then $\partial A \subseteq X_{\cl R} \subseteq M_{A(\Omega)}$.
\end{cor}

\begin{proof}  From the diagram above, we see that
the map of $\mf C(\cl R)$ into $\gamma(\mf T(\cl R))$, which is 
defined by inclusion,  yields the dotted map from $C(X_{\cl R})$ to
$C(\partial A)$.  The inclusion $\partial A \subseteq
X_{\cl R}$ follows from this. (Note  that the 
bottom row in this diagram is not exact.)
\end{proof}

Recall that we have $X_{\cl R} =  \partial A$ for $\cl R$ the 
Hardy or  Bergman module but this equality does not hold in general.
This matter  comes  down to the possible equality of $\mf C(\cl R)$
with $\gamma(\mf T(\cl R))$.

\begin{cor}\label{cor_norm}
Let $\Omega$ be a bounded domain in $\bb C^n$.
Let $\{\phi_i\}$ and $\{\psi_i\}$ be sets of functions in $A(\Omega)$,
and let $C \in \mf C(\cl R)$.  It follows that 
\[
 \big\| \sum_i M_{\phi_i} M^*_{\psi_i} + C \big\|_{\mf L(\cl R)} 
 \ge \sup_{\alpha\in\partial A } \big| \sum 
 \widehat\phi_i(\alpha)\overline{\widehat\psi_i(\alpha)} \big| .
\]
\end{cor}

Actually the inequality holds for $C$ in $\gamma(\mf T(\cl R))$.

Observe that $\mf T(\cl R)$ is the closure of the set 
\[
 \{ M_\phi M^*_\psi + \mf C(\cl R) : \phi,\psi\in A(\Omega)\} .
\]
For this conclusion we 
need show only that the set is closed under multiplication which 
follows from  the identity
\begin{align*}
 M_{\phi_1} M^*_{\psi_1} M_{\phi_2}M^*_{\psi_2} &= M_{\phi_1\phi_2} 
 M^*_{\psi_1\psi_2} + M_{\phi_1} (M^*_{\psi_1} M_{\phi_2} - M_{\phi_2} 
 M^*_{\psi_1}) M^*_{\psi_2}\\
 &= M_{\phi_1\phi_2} M^*_{\psi_1\psi_2} + M_{\phi_1} [M^*_{\psi_1}, 
 M_{\phi_2}] M^*_{\psi_2}
\end{align*}
since the second operator is in $\mf C(\cl R)$.
So this corollary is a statement about all elements of 
$\mf T(\cl R)$.

\begin{defn}
Let the \textit{Berezin spectrum} $\sigma_B(\cl R)$ of $\cl R$
denote the set of all points $\alpha$ in $M_A$ such that 
$\lim\limits_{\pmb{z}\to \alpha} \Gamma(T) (\pmb{z})$
exists for all $T \in \mf T(\cl R)$.
\end{defn}

Theorem~\ref{extend_peak} shows that $\sigma_B(\cl R)$ contains all
p-points, and thus is dense in the \v Silov
boundary $\partial A$ of $A(\Omega)$. The following result shows that
it is always a subset of the topological boundary of $\Omega$.

\begin{lem}\label{no_eig}
A quasi-free module $\cl R$ over $A(\Omega)$ has no reducing
eigenvectors.
\end{lem}

\begin{proof}
Suppose that $x$ is a reducing eigenvector for $\cl R$.
Then the functional $\Phi(A) = (Ax,x)$ is multiplicative on the
C*-algebra $\mf T(\cl R)$, and thus determines a point $\xi$ in the
character space $X_{\cl R}$.

Suppose that $\pmb{z} \in \Omega$ and $\pmb{z} \ne \xi$.
We claim that $\range V_{\pmb{z}}$ is orthogonal to $x$.
Indeed if $k = \alpha x + y$, with $y$ orthogonal to $x$, lies in
$\range V_{\pmb{z}}$, then for any $\phi \in A(\Omega)$,
\[
 \ol{\eval_{\pmb{z}}(\phi)} \alpha x +  \ol{\eval_{\pmb{z}}(\phi)}y =
 \ol{\phi(\pmb{z})} k = M_\phi^* k
 = \ol{\Phi(\phi)} \alpha x + M_\phi^* y .
\]
Moreover $M_\phi^* y$ is orthogonal to $x$, and so if $\alpha \ne
0$, $\eval_{\pmb{z}}(\phi) = \Phi(\phi)$.  
Hence $\Phi = \eval_{\pmb{z}}$ rather than evaluation at $\xi$.
This contradiction establishes the claim.

However, the span of the ranges of $V_{\pmb{z}}$ for $\pmb{z} \ne\xi$
is dense in $\cl H$.  Therefore there are no reducing eigenvalues.
\end{proof}

\begin{cor}\label{berezin_in_boundary}
The Berezin spectrum of $\cl R$ is contained in the topological
boundary of $\Omega$.
\end{cor}

\begin{proof}
The Berezin spectrum is contained in the closure of $\Omega$ in
the maximal ideal space $M_A$, which is $\Omega \cup \partial
\Omega$.  However at a point $\pmb z_0$ interior to $\Omega$, the
limit is actually the GBT evaluated at $\pmb z_0$.
For this to be multiplicative on $\mf T(\cl R)$, the range of
$V_{\pmb z_0}$ must reduce $\mf T(\cl R)$.  
The generators $M_\phi$ for $\phi$ in $A(\Omega)$ are sent to the
scalar operators $\phi(\pmb z_0) I_{\mf L(\ell^2_m)}$, and 
thus all of $\mf T(\cl R)$ is mapped into the scalars.
Therefore the range of $V_{\pmb z_0}$ must consist of reducing
eigenvalues.
This contradicts Lemma~\ref{no_eig}, and thus the
Berezin spectrum is contained in the boundary of $\Omega$.
\end{proof}

\section{Examples}\label{S:examples}

In this section, we present a variety of examples exploring the
possibilities for the character space and the GBT.

The following result shows the implication of Theorem~\ref{main_thm} in
the one  variable case in the presence  of other hypotheses on $\cl
R$. 

\begin{thm}\label{thm3}
Let $\Omega$ be a bounded domain in $\bb C$ and $\cl R$ be a 
finite rank  quasi-free Hilbert module over $A(\Omega)$. 
\begin{enumerate}

\item For $z_0$ in $\Omega$, $\range(M_{z-z_0})$ is closed if and only if 
$M_{z-z_0}$ is a Fredholm operator of index $-m$. 

\item If $\range(M_{z-z_0})$ is closed for $z_0$ in $\Omega$ 
and $m=1$ or $\mf T(\cl R)$ is irreducible, 
then $\mf K(\cl R) \subseteq \mf C(\cl R) \subseteq \gamma(\mf T(\cl R))$.

\item If $\range(M_{z-z_0})$ is closed for all $z_0$ in $\Omega$, then $\mf 
T(\cl R)/\mf K(\cl R) = C(\partial A )$, which implies 
$\mf K(\cl R) = \mf C(\cl R) = \gamma(\mf T(\cl R))$.

\end{enumerate}
\end{thm}

\begin{proof}
Statement (1) follows since the null space of any module multiplier is 
(0) and  the dimension of the null space of $(M_{z-z_0})^*$ is $m$.
Statement (2) follows  from the fact that a C*-algebra containing a
Fredholm operator of non-zero  index must contain  a non-zero compact
operator and the fact that the set of  compact operators is the minimal
ideal in any C*-algebra containing it.  Statement (3) follows from the
well-known fact that a Fredholm operator in a  C*-algebra containing
the compact operators is invertible modulo the ideal of  compact
operators.
\end{proof}

It is not always true in the above context that $\mf C = \mf K$.

\begin{eg}\label{noncompact}
Let $\ell^2$ be the Hilbert space with orthonormal basis
$\{e_k\}_{k\in\bb N}$, and $\Sigma$ be an infinite subset of 
$\bb N$ for which $\Sigma\cap (\Sigma+k)$ is finite for all $k\ge 1$ 
and $\lim\limits_{n\to \infty} \frac1n \operatorname{card}
 \{ k \in \Sigma : 0\le k < n \} = 0$.  
Let $P_\Sigma$ be the projection in $\mf L(\ell^2)$ onto the
span of $\{e_k : k\in\Sigma\}$ and $S$  be the unilateral shift on
$\ell^2$, that is, $Se_k = e_{k+1}$, $k\ge 1$. Let $A$  be 
the weighted shift $A = S\left(I-\frac12 P_\Sigma\right)$. Then the 
point  spectrum of $A^*$ is $\bb D$ and hence $A$ defines a  rank one
quasi-free  Hilbert module over $A(\bb D)$. By considering the 
polar decomposition of $A$, we see that the C*-algebra $\mf A$ 
generated by $I$  and $A$ is equal to the C*-algebra generated by $S$
and $P_\Sigma$. Moreover, $\mf A$ contains $\mf K=\mf K(\ell^2)$.

Let $\cl B = \mf A/\mf K$ and $u$ and $p$ be the images of $S$ and 
$P_\Sigma$ in $\cl B$, respectively. Then $u$ is unitary and $p$ is a
projection with the  property  that the projections
$\{u^{*k} pu^k\}_{k\in\bb Z}$ are pairwise orthogonal.  Hence,  the
ideal $\cl  J$ in $\cl B$ generated by $p$ is isomorphic to $\mf K$. 
Moreover, $\cl B$ is a trivial extension of $\mf K$ by
$C(\partial\bb D)$, that  is, $\cl B/\mf K \simeq C(\partial\bb D)$ and
the map $z\to u$ extends to an  isomorphism from $C(\partial\bb D)$ to
$\cl B$. 

Now the commutator ideal $\mf C$ of $\mf A$ contains $\mf K$ and is the 
preimage of $\cl J$ in $\mf A$. Therefore, $\mf A/\mf C\simeq
C(\partial\bb D)$ but $\mf  A/\mf K  =\cl B$ is not commutative.
Hence, $\mf K\ne \mf C$ in this case.  However, $\mf C$ and the
Berezin nullity are equal.
\end{eg}

There is another hypothesis that would imply the coincidence of 
$\gamma(\mf  T(\cl R))$ and $\mf C(\cl R)$ in the context of
Theorem~\ref{thm3}. If none of  the points in $\Omega$ lie in   
$X_{\cl R}$, then $X_{\cl R}$ must be a subset of 
$\partial\Omega$. Since $\partial\Omega \subset X_{\cl R}$, we have 
$X_{\cl R}  =\partial\Omega$ which is the maximal ideal space of   
$\mf T(\cl R)/\gamma(\mf  T(\cl R))$. 
Hence, $\mf C(\cl R) = \gamma(\mf T(\cl R))$. Note that this latter
argument is valid for the case of a domain $\Omega$ in $\bb  C^n$
if $\Omega$  is dense in $M_{A(\Omega)}$ and $\partial A =
\partial \Omega$.

 The extension of the preceding results to the several variables 
context would be of particular interest. One approach would be to 
assume that  the last stage of the Koszul complex, which defines the
Taylor spectrum, has  closed range. In this way, one can show that
$X_{\cl R} = \partial \bb B^n$ for
$\cl R$  equal to $H^2(\bb B^n)$ or $B^2(\bb B^n)$ and 
$X_{H^2(\bb D^n)} = \bb T^n$ for  all $n$. We will not provide the
details.

Now we consider when a point of the maximal ideal space lies in     
$X_{\cl R}$.  
We shall see that it can be a fairly general subset containing the
\v Silov boundary.  
Our argument works for any finitely generated algebra.
Our main application, Corollary~\ref{cor_char}, is a strengthening of
an old result of Bunce \cite{Bun} where hyponormality is assumed.
This paper was one of the precursers of Voiculescu's celebrated
Weyl-von~Neumann theorem \cite{V}  (see \cite[Lemma~II.5.5]{KD}) which
we use below.

\begin{thm}\label{character}
Let $\cl A$ be a unital operator algebra represented on a Hilbert
space $\cl H$, and $A_i, \dots, A_n$ be a finite set of generators
for $\cl A$; and let $\mf A = \textrm{C}^*(\cl A)$.
A multiplicative functional $\Phi$ on $\cl A$ extends to a character
of $\mf A$ if and only if the column operator $T$ with entries
$A_i - \Phi(A_i) I$ and $A_i^* - \ol{\Phi(A_i)}$ for $1 \le i \le n$ 
is not bounded below.
\end{thm}

\begin{proof}
Suppose that $\Phi$ extends to a character $\bar\Phi$ of $\mf A$.
If $\bar\Phi$ is nonzero on the ideal $\mf A \cap \mf K$ of compact
operators, it must correspond to a finite dimensional reducing
subspace $\cl M$ on which $T_i|_{\cl M} = \Phi(T_i) I_{\cl M}$.
Hence the column operator $T$ has $\cl M$ in its kernel.
Otherwise $\bar\Phi$ annihilates $\mf A \cap \mf K$.
So by Voiculescu's Theorem, the identity representation absorbs
countably many copies of $\bar\Phi$, that is
$\operatorname{id} \sim_a \operatorname{id} \oplus
\bar\Phi^{(\infty)}$.  This means that there is an orthonormal sequence
$\{ x_k \}$ of vectors which are asymptotically reducing and
\[
 \lim_{k\to\infty} \| A_i x_k - \Phi(A_i) x_k \| = 0 = 
 \lim_{k\to\infty} \| A_i^* x_k - \ol{\Phi(A_i)} x_k \| .
\]
So $\lim_{k\to\infty} \| T x_k \| = 0$.

Conversely, suppose that $T$ is not bounded below.
This could occur because there is a unit vector $x$ in its kernel.
Then 
\[
 A_i x = \Phi(A_i) x \qand  A_i^* x = \ol{\Phi(A_i)} x .
\]
So $x$ is a reducing eigenvalue for the operators $A_1,\dots,A_n$,
and hence for the C*-algebra $\mf A$.
Clearly restriction to $\bb C x$ is a $*$-homo\-morphism, and thus
determines a character $\bar\Phi$ of $\mf A$ which extends $\Phi$.

On the other hand, if $T$ is injective, there must be an
orthonormal sequence $\{ x_k : k \ge 1 \}$ so that 
$\lim_{k\to\infty} \| T x_k \| = 0$.
Let $P_k = x_k x_k^*$ denote the rank one projection onto $\bb C
x_k$.
Then
\begin{align*}
  \| [ A_i, P_k] \| &= \max  \big\{ \| P_k^\perp A_i P_k \|, 
  \| P_k A_i P_k^\perp \| \big\} \\ &\le  \max \big\{
 \| (A_i - \Phi(A_i)) x_k \|, \| (A_i - \Phi(A_i))^* x_k \| \big\}
\end{align*}
Hence $\lim_{k\to\infty} \| [ A_i, P_k] \| = 0$.
It follows that $\lim_{k\to\infty} \| [ A, P_k] \| = 0$ for every $A \in \mf
A$.
In particular, 
\[
 \lim_{k\to\infty} (A_i x_k, x_k ) = \Phi(A_i) \qand
 \lim_{k\to\infty} (A_i^* x_k, x_k ) = \ol{\Phi(A_i)} .
\]
Define an extension $\bar{\Phi}(A) = \lim_{k\to\infty} (A x_k, x_k )$.
It is clear that this limit exists for any polynomial in the generators
and their adjoints, and hence, by continuity, for the whole C*-algebra.
Then for $A,B \in \mf A$,
\begin{align*}
  \big|& (BA x_k, x_k) - (A x_k, x_k) (B x_k, x_k) \big| \\&=
 \big| (P_k A x_k, P_k B^* x_k) + (P_k^\perp A x_k, P_k^\perp B^* x_k)
 - (P_k A x_k, x_k) ( x_k, P_k B^* x_k) \big| 
 \\ &= \big| (P_k^\perp A x_k, P_k^\perp B^* x_k) \big|
 \le \| P_k^\perp A P_k \| \, \| P_k^\perp B^* P_k \|  .
\end{align*}
Hence $\bar{\Phi}(BA) = \bar{\Phi}(A) \bar{\Phi}(B)$ is multiplicative.
So $\bar{\Phi}$ is a character of $\mf A$ which extends $\Phi$.
\end{proof}

By the previous theorem and Lemma~\ref{no_eig}, we obtain:

\begin{cor}\label{cor_char}
Let $\cl R$ be a quasi-free module over $A(\Omega)$ where $\Omega$ is
a domain in $\bb C^n$.  Then evaluation at a point $\omega \in
\Omega$ extends to a multiplicative functional on $\mf T(\cl R)$ if
and only if the column operator $T$ with entries
$M_{z_i - \omega_i}$ and $M_{z_i - \omega_i}^*$ for $1 \le i \le
n$ is not bounded below.
\end{cor}

\begin{eg}\label{wtd_shift}
Consider a contractive rank one module over $A(\bb D)$ which is rotation
invariant in the sense that sending $z$ to $\lambda z$, where
$|\lambda|=1$, induces an isometric automorphism of $\cl R$.  It is
routine to verify that this implies that $\{ z^k : k \ge 0 \}$ is an
orthogonal set spanning $\cl R$. Hence $M_z$ is an injective weighted
shift.  Being a representation of
$A(\bb D)$, we see that $\| M_z \| = 1 = \spr(M_z)$ (where $\spr(T)$
is the spectral radius of $T$).
All weighted shifts are unitarily equivalent to a shift with positive
weights, so we may suppose that $M_z \simeq T$ where $T$ acts on $\ell^2$
by $Te_n = a_n e_n$ and $0 < a_n \le 1$ for $n \ge 0$.

The problem of determining the set $X_{\cl R}$ for this module comes down
to determining when the column operator with entries $T-\lambda I$ and
$(T-\lambda I)^*$ is bounded below.  Equivalently, we may consider the
operator 
\[
 X = (T-\lambda I)^*(T-\lambda I) + (T-\lambda I)(T-\lambda I)^* .
\]
Suppose that for each $\ep>0$ and $d \in \bb N$, there is an integer $k$ so
that 
\[
 \big| a_{k+i} - |\lambda| \big| < \ep \qfor  0 \le i \le d+1 ,
\]
then with $\lambda = e^{i\theta} |\lambda|$, we may set
$x_k = d^{-1/2} \sum_{j=1}^d e^{-ij\theta} e_{k+i}$.
It is a simple exercise to show that $\| (T-\lambda I) x_k \|$ and
$\| (T-\lambda I)^* x_k \|$ are both small.  In the limit, we see that the
column operator is not bounded below and $\lambda$ lies in $X_{\cl
R}$.

On the other hand, if $a_k$ is bounded away from $|\lambda|$ for all $k
\ge0$, say by $\delta>0$, then $X$ is invertible.
Indeed, one can easily see that 
\[
 X = T^*T + TT^* + 2 |\lambda|^2 I - 2 (\ol{\lambda}T + \lambda T^*) .
\]
This is a tridiagonal operator with entries
\[
 x_{ii} = a_{i-1}^2 + a_i^2 + 2 |\lambda|^2 \qquad
 x_{i,i+1} = 2 \ol{\lambda} a_i \qand 
 x_{i+1,i} = 2 \lambda a_i .
\]
(We set $a_{-1}=0$.)  Then $X$ is the sum of two operators, 
$X_0$:
\[
 \begin{bmatrix}
  a_0^2 \!+\! |\lambda|^2 & 2 \ol{\lambda} a_0 &&&&&\\
  2 \lambda a_0 & a_0^2 \!+\! |\lambda|^2 &&&&&\\
  && a_2^2 \!+\! |\lambda|^2 & 2 \ol{\lambda} a_2 &&&\\
  && 2 \lambda a_2 & a_2^2 \!+\! |\lambda|^2 &&&\\
  &&&& a_4^2 \!+\! |\lambda|^2 & 2 \ol{\lambda} a_4 &\\
  &&&& 2 \lambda a_4 & a_4^2 \!+\! |\lambda|^2 &\\
  &&&&&& \ddots
 \end{bmatrix}
\]
and $X_1$:
\[
 \begin{bmatrix}
  |\lambda|^2 &&&&&&&\\
  &a_1^2 \!+\! |\lambda|^2 & 2 \ol{\lambda} a_1 &&&&&\\
  &2 \lambda a_1 & a_1^2 \!+\! |\lambda|^2 &&&&&\\
  &&& a_3^2 \!+\! |\lambda|^2 & 2 \ol{\lambda} a_3 &&&\\
  &&& 2 \lambda a_3 & a_3^2 \!+\! |\lambda|^2 &&&\\
  &&&&& a_5^2 \!+\! |\lambda|^2 & 2 \ol{\lambda} a_5 &\\
  &&&&& 2 \lambda a_5 & a_5^2 \!+\! |\lambda|^2 &\\
  &&&&&&& \ddots
 \end{bmatrix}
\]
It is now easy to check that $X \ge X_0 \ge \delta^2/2 \, I$.
So $\lambda$ does not lie in $X_{\cl R}$.

A similar analysis, which we omit, shows that if there is a $\delta >0$
and an integer $d$ so that there are no strings of weights of length $d$
all within $\delta$ of $|\lambda|$, then again $X$ is bounded below;
and so $\lambda$ does not lie in $X_{\cl R}$.
Various spectral properties of weighted shifts have been studied in detail.
See Shields \cite{Sh} for an overview.  This result could well be known, 
but we have been unable to find it in the literature.

Let $C$ be a compact subset of $[0,1]$ containing $1$ in which $0$ is not
an isolated point.  Choose a sequence $\{ r_k : k \ge 1 \}$ contained in
$C\bsl\{0\}$ so that every point in $C$ is a cluster point.
Then define $T$ to be the weighted shift with weights 
$r_1,r_2,r_2,r_3,r_3,r_3,r_4,r_4,r_4,r_4,\dots$.
In this case, the analysis above shows that $X_{\cl R} = C \bb T$.
This  set of examples shows that $X_{\cl R}$ can be a rather arbitrary
subset of $M_A$ containing $\partial M_A$.

For another example, if $0<r<1$, let the weights $a_n$ slowly oscillate
back and forth between $1$ and $r$.  That is, $|a_i - a_{i+1}|$ tends to
$0$.  Then $T$ is essentially normal with essential spectrum equal to the
annulus $\{ \lambda : r \le |\lambda| \le 1 \}$.  For this example, it is
clear that every point in the annulus is an approximate reducing
eigenvalue because $T$ is approximately unitarily equivalent to $S \oplus
N$ where $S$ is the unilateral shift and $N$ is normal and has spectrum
equal to the annulus.
\end{eg}

In the following example, we show that the disk algebra has a
completely bounded module $\cl R$ (hence similar to a contractive one)
in which the set $X_{\cl R}$ is empty.  Indeed, it will be similar
to the classical Hardy module.  This should be compared to 
Corollary~\ref{Shilov} which implies that a quasi-free
(hence contractive) module over the unit disk must contain the
whole unit circle in $X_{\cl R}$.

\begin{eg}\label{simple}
Fix $0 < r < 1$ and consider the weighted shift $T$ with weights
$1,r^{1/2},1,r^{1/4},1,r^{1/2},1,r^{1/8},\dots$, namely
$a_n = r^{1-\operatorname{gcd}(n+1,2^n)^{-1}}$.
Note that $T^{2^k}$ shifts by $2^k$ with weights $b_n$ which are the
product of $2^k$ successive weights of $T$, of which 
$2^{k-1}$ are 1, $2^{k-2}$ are $r^{1/2}$, \dots, one is
$r^{1-2^{1-k}}$ and one lies between $r$ and $r^{1-2^{-k}}$.
So there is an integer $p = p(i) \ge n$ so that
\begin{align*}
 2^{-k} \log_r b_i &= 1 \!-\! \Big( \tfrac12(1) \!+\! 
\tfrac14\big(\tfrac12\big) 
 \!+\! \tfrac18 \big(\tfrac34\big) \!+\! \dots \!+\!
 2^{-k}\big(1\!-\!2^{1-k}\big) \!+\! 2^{-k}\big(1\!-\!2^{-p}\big) \Big)
 \\&= \frac13 \big( 1 \!+\! 2^{-k-p} \!-\!\tfrac23 4^{-k} \big) .
\end{align*}
Therefore these weights satisfy
\[
    r^{2^{1-k}/3} < r^{-2^k/3} b_n \le r^{-2^{-k}/3} .
\]
In particular, the spectral radius of $T$ is $r^{1/3}$.

Define $A = r^{-1/3} T$.  The previous estimate establishes that  
$$\| A^{2^k} \| \le r^{-2^{-k}/3}.$$
Hence it is a routine calculation to see that $A$ is power bounded.
Likewise, all powers of $A$ are uniformly bounded below.
By \cite{Sh}, $A$ is similar to the unilateral shift.
Thus $A$ determines a bounded module $\cl R$ which is similar to the
Hardy module.
Nevertheless, we shall show that $X_{\cl R}$ is empty.

It is evident that $\mathrm{C}^*(A) = \mathrm{C}^*(T)$ is
generated by the unilateral shift $S$ and the diagonal operator $D =
(T^*T)^{1/2}$, both of which come from the polar decomposition of $T$.
The spectrum of $D$ is
\[
 \sigma(D) = \{ r, r^{1-\operatorname{gcd}(n+1,2^n)^{-1}}
 : n \ge 0 \} ;
\]
 and the spectral projection for the point 
$r^{1-\operatorname{gcd}(n+1,2^n)^{-1}}$ is the projection onto  
\[
 \{ e_i : i+1 \equiv 2^k \,(\operatorname{mod} 2^{k+1}) \} .
\]
It is clear that conjugating these projections by powers of $S$
yields the projections $E_{k,l}$ onto the subspaces spanned by
\[
 \{ e_i : i \equiv l \,(\operatorname{mod} 2^k) \}
 \quad \text{ for } k \ge 1 \text{ and } 0 \le l < 2^k .
\]
It is easy to see that this is the C*-algebra generated by all
weighted shifts of period $2^k$ for all $k \ge0$.
This is an extension of the compact operators by the $2^\infty$
Bunce--Deddens algebra \cite[\S V.3]{KD}.

The commutator ideal properly contains the compact operators because
$A$ is not essentially normal.  As the Bunce--Deddens algebras are
simple, the commutator ideal is the whole C*-algebra.
In particular, there are no characters.
Also, by the analysis of Example~\ref{wtd_shift}, we can see by a
different method that there are no characters.
\end{eg}

We can also extend our results to cover spherical contractions. 
Recall that the Hilbert module $\cl M$ over $A(\bb B^n)$ is said to be
{\em  spherically contractive\/} if 
$\sum\limits^n_{i=1} \|M_{z_i}f\|^2_{\cl R} \le \|f\|_{\cl R}$ for all
$f$ in $\cl R$. 

\begin{thm}\label{spherical}
If $\cl R$ is a finite rank quasi-free Hilbert  module over $A(\bb B^n)$ which 
is spherically contractive, then the GBT defines a $*$-homomorphism 
from $\mf  T(\cl R)$ onto $C(\partial\bb B^n)$ with kernel ideal
$\gamma(\mf T(\cl R))$  containing $\mf C(\cl R)$.
\end{thm}

\begin{proof}
The first part of the proof is similar to that of Lemma~\ref{asymp_reducing}
but the  operator $M_\phi$ on $\cl R$ is replaced by a column of such
operators.  Let  $\Phi = \begin{spmatrix}\phi_1\\ \vdots\\ \phi_k
\end{spmatrix}$  be a $(1\times k)$-column of functions in $A(\bb
B^n)$,  and let $M_\Phi$ be the corresponding $(1\times k)$-column
matrix of operators  acting from $\cl  R$ to $\cl R^{(k)}$, the direct
sum of $k$ copies of $\cl R$. An easy  calculation shows that
\[
 \Big\|\sum_i M^*_{\phi_i} M_{\phi_i}\Big\|_{\mf L(\cl R)} = \|M^*_\Phi 
 M_\Phi\|_{\mf L(\cl R)} = \|M_\Phi\|^2 
\]
and
\[
 \sup_{\pmb{z}\in \bb B^n} \|\Phi(z)\|^2_{\mf L(\ell^2_m,\ell^{2\,(k)}_m)}
 = \sup_{\pmb{z}\in \bb B^n} \sum_i |\phi_i(\pmb{z})|^2.
\]
Now for the column matrix operator $M_Z$ from $\cl R$ to $\cl R^{(n)}$, defined 
by the column matrix 
$Z = \begin{spmatrix} z_1\\ \vdots\\ z_n\end{spmatrix}$ formed
from the coordinate functions, we have 
\[
 \| M_Z \| \le \sup_{\pmb{z}\in \bb B^n} \sum_{i=1}^n |z_i|^2 = 1
\]
and hence $M_Z$ is a spherical contraction. 

Next decompose $\cl R = \cl P_{\pmb{z}} \oplus \cl 
Q_{\pmb{z}}$ for $\pmb{z}$ in $\bb B^n$, 
where $\cl P_{\pmb{z}} = \range  V_{\pmb{z}}$ and $\cl Q_{\pmb{z}}$ is the
orthogonal complement of $\cl P_{\pmb{z}}$. 
Similarly, $\cl R^{(n)} = \cl P^{(n)}_{\pmb{z}} \oplus \cl 
Q^{(n)}_{\pmb{z}}$. 
Moreover, if we decompose the matrix for the operator $M_Z$ relative 
to this  decomposition, we obtain $M_Z \simeq \begin{sbmatrix}
A_{\pmb{z}}&0\\  B_{\pmb{z}}&C_{\pmb{z}}\end{sbmatrix}$, where
$A_{\pmb{z}}$ is the  column matrix $\begin{spmatrix}  z_1I_{\cl
P_{\pmb{z}}}\\ 
\vdots\\ z_n I_{\cl P_{\pmb{z}}}\end{spmatrix}$. Since $M_Z$ is a 
contraction, by the theorem of Sz.-Nagy and Foia\c{s}, we have that 
$B_{\pmb{z}} =  (I-A^*_{\pmb{z}} A_{\pmb{z}})^{1/2} B'_{\pmb{z}}$,
where $B'_{\pmb{z}}$ is a  contraction. 
But 
\[
 (I-A^*_{\pmb{z}}A_{\pmb{z}})^{1/2} = (1-\sum |z_i|^2)^{1/2}
\] 
which implies that 
${\overline{\lim}}_{\|\pmb{z}\|\to 1} \|B_{\pmb{z}}\|  = 0$.

This result enables us to establish the conclusion of
Proposition~\ref{extend_peak} for functions in the
closed linear span of $\{z_i\}$ which  coincides with $A(\bb B^n)$ and
this concludes the proof.
\end{proof}

An example of a Hilbert module for which Theorem~\ref{spherical} applies
but not  Theorem~\ref{main_thm} 
is the non-commutative Hardy space $H^2_d$
which is a module  over $\bb C[\pmb{z}]$ for which the coordinate
functions are contractive but 
$H^2_d$ is not a bounded module over $A(\bb B^n)$.

\section{Extension to a Larger Algebra}\label{sec4}

On the disk $\bb D$ one knows that $H^2(\bb D)$ and $B^2(\bb D)$ are 
modules not  only for $A(\bb D)$ but for the larger algebra,
$H^\infty(\bb D)$,  of all  bounded holomorphic functions over $\bb
D$. The same thing is true for other  quasi-free Hilbert modules. In
this section we restrict our attention to the  rank 
one case, $m=1$, although most of the following developments carry 
over to the  higher rank case. Many of the techniques are well-known
and often are part of  the folklore.

Let $\cl R$ be a rank one quasi-free Hilbert module over $A(\Omega)$, 
$\Hol(\Omega)$ be the algebra of all holomorphic functions on $\Omega$, 
and $H^\infty(\Omega)$ be the subalgebra of bounded functions in 
$\Hol(\Omega)$. A function $\phi$ in $\Hol(\Omega)$ is said to be 
a multiplier for $\cl R$ if $\phi\cl R\subset \cl R$. Let 
$\Mul(\Omega)$ denote the algebra of all multipliers for $\cl R$ with norm 
$\|\phi\|_{\Mul(\cl R)} = \|M_\phi\|$, where $M_\phi$ is the 
operator on $\cl R$ defined to be multiplication by $\phi$, which is in $\mf 
L(\cl R)$ by the closed graph theorem. Standard arguments show that
not only must multipliers be bounded, but that
$\|M_\phi\|_{\mf L(\cl R)} \le  \|\phi\|_{H^\infty(\Omega)}$.
Hence, $\Mul(\cl R)$ is a subalgebra of  $H^\infty(\Omega)$ which
is contractively included. In general, the norms are  not equal
and $\Mul(\Omega) \ne H^\infty(\Omega)$. The following lemma 
shows that $\Mul(\cl R)$ is closed in a rather weak topology.

A sequence $\{\phi_i \}$ in $A(\Omega)$ is said to converge in the 
bounded,  pointwise limit (bpwl) topology to a function $\phi$ in
$H^\infty(\Omega)$ if
\begin{enumerate}
 \item $\sup_i \|\phi_i\|_{A(\Omega)} <\infty$ and
 \item $\lim_i \phi_i(\pmb{z}) = \phi(\pmb{z})$ for $\pmb{z}$ in 
$\Omega$.
\end{enumerate}

\begin{lem}\label{lem4.1}
Let $\cl R$ be a quasi-free Hilbert module over $A(\Omega)$.
Suppose that $\phi_i \in \Mul(\cl R)$ converge bpwl to $\phi$.
Then $M_{\phi_i}$ converge in the weak operator topology to $M_\phi$.
Thus $\Mul(\cl R)$ is closed in the bpwl-topology on $H^\infty(\Omega)$.
Moreover $\{ M_\phi : \phi \in \Mul(\cl R), \|\phi\|_{H^\infty(\Omega)} \le 1
\}$ is weak operator closed.  
\end{lem}

\begin{proof}
Let $\phi$ be in $H^\infty(\Omega)$ and $\{\phi_i\}$ be a sequence in 
$\Mul(\cl R)$ which converges to $\phi$ in the bpwl-topology. 
Since the  sequence of operators $\{M_{\phi_i}\}$ is uniformly bounded, a
subsequence of it converges to some $X$ in $\mf L(\cl R)$ in the weak 
operator topology.
For a  vector $x$ in $\ell^2_m$ and $g$ in $\cl R$, we have
\begin{align*}
\langle\eval_{\pmb{z}}(Xg),x\rangle_{\ell^2_m} &= \langle Xg, 
V_{\pmb{z}}x\rangle_{\cl R} = \lim_i \langle 
M_{\phi_i}g,V_{\pmb{z}}x\rangle_{\cl R}\\
&= \lim_i \langle g, \overline{\phi_i(\pmb{z})} V_{\pmb{z}} x\rangle = 
\langle g,\overline{\phi(z)} V_{\pmb{z}} x\rangle_{\cl R} 
\\&= \langle 
\eval_{\pmb{z}}(\phi g), x\rangle_{\ell^2_m}.
\end{align*}
Therefore $Xg = \phi g$; whence $X = M_\phi$ and $\phi$ is in 
$\Mul(\cl R)$. Since this was valid for \emph{any} subsequence which
had a limit in the weak operator topology, standard arguments show
that the original sequence converge to $M_\phi$ in the weak operator
topology.

Conversely, suppose that $\phi_i \in \Mul(\cl R)$ such that 
$M_{\phi_i}$ converge in the weak operator topology to $X$.
Then since $\Gamma_{\pmb{z}}(X) = V^*_{\pmb{z}} X V_{\pmb{z}}$ is 
compressed on both sides to the finite dimensional range of
$V_{\pmb{z}}$, it follows that
\[
 \Gamma_{\pmb{z}}(X) = \lim_{i \to \infty} \Gamma_{\pmb{z}}(M_{\phi_i})  
  = \lim_{i \to \infty} \phi_i(\pmb{z}) I_n.
\]
Therefore there is a bounded pointwise limit $\phi$ of the sequence 
$\phi_i$. So $\phi \in H^\infty(\Omega)$ and is a bpwl-limit of
multipliers, whence lies in $\Mul(\cl R)$ by the previous paragraph.
Repeating the computation above, it follows that $X = M_\phi$.
\end{proof}

An argument along the same lines is given in \cite[Thm.~5.2]{C-S}.

If the inclusion is isometric, namely 
$\| \phi \|_{\Mul(\cl R)} = \| \phi\|_{H^\infty(\Omega)}$, then 
$\{ M_\phi : \phi \in \Mul(\cl R), \|\phi\|_{H^\infty(\Omega)} \le 1 \}$
is the unit ball of the algebra of multipliers
$\{ M_\phi : \phi \in \Mul(\cl R) ) \}$.
In this case, the Krein--Smulyan Theorem shows that this space is 
weak-$*$ closed.

\begin{pro}\label{pro4.3}
Let $\Omega$ be a domain in $\bb C^n$ for which $A(\Omega)$ is dense in 
$H^\infty(\Omega)$ in the bpwl-topology. Then every quasi-free Hilbert module 
over $A(\Omega)$ extends to a bounded Hilbert module over 
$H^\infty(\Omega)$. 
\end{pro}

\begin{proof}
By Lemma~\ref{lem4.1} and the hypotheses we see that $H^\infty(\Omega) = 
\Mul(\Omega)$. Therefore, we can define a module action 
$H^\infty(\Omega) \times \cl R\to\cl R$. By the closed graph theorem,
there exists a constant $K>0$ such  that 
$\|M_\phi\|_{\mf L(\cl R)}\le  K\|\phi\|_{H^\infty(\Omega)}$ for 
$\phi$ in $H^\infty(\Omega)$.
\end{proof}

It is known that various assumptions on $\Omega$ imply that $A(\Omega)$ 
is dense  in $H^\infty(\Omega)$ in the bpwl-topology. We provide one
easy example. Recall  that a domain $\Omega$ in $\bb C^n$ is said to
be {\em starlike relative to a  point\/} $\pmb{z}_0$ in $\Omega$, if
for every point $\pmb{z}$ in $\Omega$ the  line segment with endpoints
$\pmb{z}_0$ and $\pmb{z}$ is contained in $\Omega$.  Moreover,
$\Omega$ is said to be {\em starlike\/} if it is starlike relative to
some point  in $\Omega$.

If we want $\cl R$ to be a contractive Hilbert module over 
$H^\infty(\Omega)$,  we need more information about how $A(\Omega)$ is
bpwl-dense in $H^\infty(\Omega)$. The most straightforward hypothesis
is the  assumption that the sequence $\{\phi_i\}$ in $A(\Omega)$
converging in the  bpwl-topology to a function $\phi$ in
$H^\infty(\Omega)$ can be chosen such  that 
$\|\phi_i\|_{A(\Omega)} \le \|\phi\|_{H^\infty(\Omega)}$. 
We will say  that $A(\Omega)$ is {\em
strongly bpwl-dense\/} in  $H^\infty(\Omega)$ in that case.

\begin{lem}\label{lem4.2}
If $\Omega$ is a starlike domain, then $A(\Omega)$ is strongly 
bpwl-dense in $H^\infty(\Omega)$.
\end{lem}

\begin{proof}
Let $\pmb{z}_0$ be a point in $\Omega$ relative to which it is starlike. With 
$\pmb{z}_0$ as origin and $0\le t\le 1$, let $F_t$ be the $t$ dilation 
of $\bb  C^n$; that is, $F_t(\pmb{z}) = \pmb{z}_0 +
t(\pmb{z}-\pmb{z}_0)$ for $\pmb{z}$  in 
$\bb C^n$. By the starlike hypothesis, we have 
$F_t(\Omega)\subset \Omega$ for $0<t\le 1$. 

For $\phi$ in $H^\infty(\Omega)$, if $\phi_t = \phi\circ F_t$, then 
$\{\phi_{1-\frac1k}\}$ is a sequence in $A(\Omega)$ for which 
$\|\phi_{1-\frac1k}\|_{H^\infty(\Omega)} \le \|\phi\|_{H^\infty(\Omega)}$ and 
$\lim\limits_{k\to\infty} \phi_k(\pmb{z}) = \phi(\pmb{z})$ for $\pmb{z}$ 
in $\Omega$. Thus $A(\Omega)$ is dense in $H^\infty(\Omega)$ in the
bpwl-topology.
\end{proof}

The arguments given in Section~\ref{sec3} carry over to the algebra of 
multipliers and hence to $H^\infty(\Omega)$ in case $A(\Omega)$ is 
strongly  dense in the bpwl-topology. We will  state the results in
one case and comment only on any  necessary 
changes.

\begin{defn}\label{defn4.4}
A domain $\Omega$ for which $A(\Omega)$ is  bpwl-dense in 
$H^\infty(\Omega)$ is  said to be {\em weakly pointed\/} at $\alpha$
in the \v Silov boundary of $M_{H^\infty(\Omega)}$ if
$H^\infty(\Omega)$ is pointed at $\alpha$.
\end{defn}

Again the Izzo--Wermer result applies to show that $H^\infty(\Omega)$
is pointed at all p-points, which forms a dense subset of the \v Silov
boundary.

Let $\Omega$ be a bounded domain in $\bb C^n$ for which $A(\Omega)$ is 
strongly bpwl-dense in $H^\infty(\Omega)$.
If $\cl R$ is a finite rank  quasi-free Hilbert module over
$A(\Omega)$, then $\cl R$ is a contractive  Hilbert module over
$H^\infty(\Omega)$. Let $\mf T_\infty(\cl R)$ be the  C*-algebra
generated by $\{M_\phi : \phi\in H^\infty(\Omega)\}$, $\mf 
C_\infty(\cl R)$ be the commutator ideal in $\mf T_\infty(\cl R)$.
Then $\mf  T_\infty(\cl R)/\mf C_\infty(\cl R) \simeq
C(X^\infty_{\cl R})$ for a subset
$X^\infty_{\cl R}$ of $M_{H^\infty(\Omega)}$. Finally, let 
$\Gamma_\infty$  be the restriction of the GBT to $\mf T_\infty(\cl
R)$ and $\sigma_\infty$ the $*$-homomorphism from $\mf T_\infty(\cl
R)$ to $C(X^\infty_{\cl R})$.

\begin{thm}\label{thm8}
Let $\Omega$ be a bounded domain in $\bb C^n$ for which $A(\Omega)$ is 
strongly  bpwl-dense in $H^\infty(\Omega)$. If $\cl R$ is a  finite 
rank quasi-free Hilbert module over $A(\Omega)$, then the GBT  defines
a $*$-homomorphism from $\mf T_\infty(\cl R)$  onto $\mf C(\partial
H^\infty(\Omega))$ whose kernel contains $\mf  C_\infty(\cl R)$.
\end{thm}

\begin{proof}
The proof is the same  as Theorem~\ref{main_thm}
with $H^\infty(\Omega)$ replacing $A(\Omega)$.
The same problem exists in this case regarding the possible
continuity of $\Gamma_\infty$ on $\Omega \cup \partial
H^\infty(\Omega)$.
\end{proof}

The classic examples of the foregoing structure occur for $\Omega$
the  unit  disk and $\cl R$ the Hardy and Bergman modules. In
\cite{D}, the second author  established that the quotient 
$\mf T_\infty(H^2(\bb D))/\mf C_\infty(H^2(\bb  D))$ is 
isometrically isomorphic to 
$L^\infty(\bb T)\simeq C(\partial H^\infty(\bb D))$ 
via the symbol map which coincides on $\partial H^\infty(\bb D)$ with 
the GBT.  
Hence, in this case the extension of the GBT restricted to
$\partial H^\infty(\bb D)$ coincides with the symbol map and the
Berezin nullity coincides with $\mf C_\infty(H^2(\bb D))$. cf.
\cite{A-Z} The result of McDonald and Sundberg \cite{M-S} shows
that 
$\mf T_\infty(B^2(\bb D))/\mf C_\infty(B^2(\bb  D))$  is isometrically
isomorphic to $C(M_1)$, where $M_1$ is the subset of $M_{H^\infty(\bb D)}$
consisting of the one-point parts. Since $\partial  H^\infty(\bb
D)\subsetneq M_1$, we see in this case that the kernel of the
extension of the GBT restricted to $\partial H^\infty(\bb D)$ does
not equal $\mf C_\infty(B^2(\bb D)))$.

There is a natural family of kernel Hilbert spaces $\cl R_n$ on $\bb D$ 
which  define quasi-free Hilbert modules over $A(\bb D)$. The kernel
functions are 
$(1-z\overline\omega)^{-(n+1)}$ and the first two are  the Hardy 
and Bergman modules. A natural problem is  to determine the maximal 
ideal space of $\mf T_\infty(\cl R_n)/\mf C_\infty(\cl R_n)$ which can 
be  identified with a closed subset $M_n$ of $M_{H^\infty(\bb D)}$.
Since these kernels are  each invariant under the conformal self-maps
of $\bb D$, it follows that the  maximal ideal spaces $M_n$ are also.
Another interesting question would seem to be the  characterization of
all  conformally invariant closed subsets of $M_{H^\infty(\bb D)}$, or
at  least the ones that are the maximal ideal space of  subalgebras
between $H^\infty(\bb D)$ and $L^\infty(\bb T)$ or, equivalently, the
conformally  invariant Douglas algebras. Known examples of the latter
are $H^\infty(\bb D)$, $H^\infty(\bb D)+C$, $A_{\text{int}}$  (the
algebra  generated by $H^\infty(\bb D)$ and the complex conjugate of
the interpolating  Blaschke products) and $L^\infty(\bb T)$, with the
corresponding maximal ideal spaces being $M_{H^\infty(\bb D)}$,
$M_{H^\infty(\bb D)}\setminus \bb D$, $M_1$, $\partial  H^\infty(\bb
D)$. An example not on this list is obtained from the Douglas  algebra
generated by $H^\infty(\bb D)$ and the complex conjugates of the 
singular inner functions for purely atomic measures. It would be of
interest to  better understand this algebra.

\section{Multiplication Operators Bounded Below}\label{sec5}

In this final section we want to relate further the commutator ideal 
and the  Berezin nullity to the nature of the multiplication
operators. We have already considered one aspect of such a
relationship in Theorem~\ref{thm3} and the paragraphs that  follow.
There we considered  the implication of the operators $M_{z-z_0}$
having closed range for  the case of $\Omega$ a planar domain. We want
to consider this matter in more detail and restrict  attention, at
first, to the case of the unit disk. We first recall that  the 
question can be reduced to that of inner functions.

\begin{lem}\label{lem5.1}
Let $\cl R$ be a finite-rank quasi-free Hilbert module over $A(\bb D)$, 
$\phi$ in $H^\infty(\bb D)$, and $\phi=\theta f$, where $\theta$ is 
inner  and $f$ is outer. The operator $M_\phi$ has closed range if and
only if $f$ is invertible in $H^\infty(\bb D)$ and $M_\theta$ has
closed range.
\end{lem}

\begin{proof}
First observe that if $f$ is an outer function, then $M_f$ has
dense range.  
In \cite[Theorem~II.7.4]{G} it is shown that if $f$ is outer, there are
functions  $g_n \in H^\infty(\bb D)$ so that $\|fg_n\|_\infty \le 1$ and
$fg_n$ converges to $1$ almost everywhere on the unit circle.
Since point evaluation in the disk is absolutely continuous, and the
kernels are continuous on $\bb D$, it follows that $f(z) g_n(z)$
converges pointwise to $1$; and in fact it converges uniformly on
compact subsets of $\bb D$.
In particular it converges in the bpwl topology, and so it follows from
Lemma~\ref{lem4.1} that $M_{fg_n}$ converges to the identity in the weak
operator topology.
As the range of $M_f$ contains the range of $M_{fg_n}$ for all $n \ge 1$,
one deduces that $M_f$ has dense range.

Since $M_\phi = M_\theta M_f$ has closed range, the operator is bounded 
below because it is one-to-one. 
Thus $M_f$ is bounded below.  
Since it has dense range, $M_f$ is invertible;
and so $M_\theta$ has closed range.

Consequently $M_f^*$ is invertible, and hence bounded below.
So for every vector $k_z$ in the range of $V_z$,
\[
 \| \overline{f(z)} k_z\| = \| M_f^* k_z \| \ge \ep \| k_z \| .
\]
Whence $|f(z)| \ge \ep$ for $z \in \bb D$.
Therefore,  $1/f$ is in $H^\infty(\bb D)$.
\end{proof}

\begin{pro}\label{pro5.2}
Let $\cl R$ be a finite-rank quasi-free Hilbert module over $A(\bb D)$ and 
let $\phi$ be in $H^\infty(\bb D)$. If $M_\phi$ has closed range,
then $\widehat\phi$ is non-zero on $X^\infty_{\cl R}$.
\end{pro}

\begin{proof}
Since $M_\phi$ is one-to-one, $M_\phi$ has closed range if 
and  only if $M^*_\phi M_\phi$ is an invertible operator. In that
case, it is  invertible in $\mf T_\infty(\cl R)$ and hence its
image, $\hat{\phi}|_{X^\infty_{\cl R}}$, is invertible in
$C(X^\infty_{\cl R})$ which  completes the proof.
\end{proof}

An interesting question is whether the converse of this result holds, 
that is,  whether the maximal ideal space of $\mf T_\infty(\cl R)/\mf
C_\infty(\cl R)$  is determined by where the Gelfand transforms of
the  functions in $H^\infty(\bb  D)$ whose 
multipliers have closed range, are non-zero. First, the converse 
holds for $\cl R = H^2(\bb D)$. The argument proceeds as follows:\ 
Multiplication by each inner function defines an isometry on $H^2(\bb
D)$ and hence has closed  range. Finally, the \v Silov boundary of
$H^\infty(\bb D)$ consists precisely of  those points 
in $M_{H^\infty(\bb D)}$ for which the Gelfand transforms of all inner functions 
don't vanish.  In 
the case $\cl R = B^2(\bb D)$, Horowitz \cite{H} and McDonald and 
Sundberg \cite{M-S} have shown that the same thing holds, that is,
$X^\infty_{B^2(\bb D)}$ is precisely the subset of 
$M_{H^\infty}(\bb D)$ on which an inner function is non-zero precisely
when as a module multiplier it has  closed range. Moreover, they showed
that the inner functions with closed range on $B^2(\bb D)$ are 
precisely the  finite products of  interpolating Blaschke products.

Also, one might hope that closedness of the range for $M_\phi$ could 
be  inferred directly from the Gelfand transform (or the Berezin
transform) of the function $\phi$. The preceding analysis shows that
is the case for $\cl R = B^2(\bb  D)$ if one considers the asymptotic
behavior of the transform on $M_{H^\infty(\bb D)}$ in the 
directions of $X^\infty_{B^2(\bb D)}$ or at limit nets that converge
to the one-point parts of $M_{H^\infty(\bb  D)}$. There
is another more direct possibility, however. If $M_\phi$ has closed 
range,  then $\|M_\phi k_z\| \ge \ep\|k_z\|$ for some $\ep>0$ and all 
$z$ in $\bb D$, where $k_z$ is in the range of $V_z$. Is the converse 
true?  That is, is it enough to test whether $M_\phi$ is bounded below
on just these  vectors? This is trivially the case for $\cl R =
H^2(\bb D)$. Zhu has informed  us that this is the case for $\cl R =
B^2(\bb D)$ for $\phi$ an inner function. His proof rests on the
characterization of the  inner functions on $B^2(\bb D)$ with closed
range. What about for other  quasi-free Hilbert modules over $A(\bb
D)$?

If one examines the proof of the preceding proposition, one sees that it extends 
to a more general setting with no change.

\begin{pro}\label{pro5.3}
Let $\Omega$ be a bounded domain in $\bb C^n$. 
For $\cl  R$ a finite-rank quasi-free Hilbert module over $A(\Omega)$
and $\phi$ in $H^\infty(\Omega)$, a necessary condition for $M_\phi$
to have closed range  is for $\hat{\phi}|_{X^\infty_{\cl R}}$ to be
invertible in $C(X^\infty_{\cl R})$.
\end{pro}

In \cite{S} Sundberg obtains a converse to this result under somewhat 
different  but related assumptions. In particular, he assumes that
$\cl R$ is the compression  of a larger module $N$ on which the module
action of $A(\Omega)$ is via  normal operators. We adapt his proof to
our context but at the price of having  to consider finite columns of
operators defined on $\cl R$ by module  multiplication or their
adjoints.

Based on Theorem~\ref{thm8} we can assume $\sigma_\infty$ is defined 
on the  whole \v Silov boundary of $M_{H^\infty(\Omega)}$.
The following result is an extension of Theorem~\ref{character} to the
case when there are not finitely many generators.

\begin{thm}\label{thm9}
Let $\Omega$ be a bounded domain in $\bb C^n$, 
and let $\cl R$ be a finite-rank quasi-free Hilbert module over
$A(\Omega)$. A point $\alpha$ in $M_{H^\infty(\Omega)}$ is in
$X^\infty_{\cl R}$ if and only if  all operators of the form 
$T = \begin{sbmatrix}T_1\\ \vdots \\ T_n\end{sbmatrix}$
mapping $\cl  R$ to $\cl R^{(n)}$ with $T_i = M_{\phi_i}$ or 
$M^*_{\phi_i}$ for $\phi_i$ in $H^\infty(\Omega)$ satisfying
$\hat{\phi}_i(\alpha) = 0$ have non-closed range. 
\end{thm}

\begin{proof}
Fix $\alpha$ in $X^\infty_{\cl R}$, and let $\Phi$ be the
character on $\mf T_\infty( \cl R)$ extending evaluation at
$\alpha$.  Consider an operator 
$T  = \begin{bmatrix}T_1\\ \vdots\\ T_k \end{bmatrix}$ mapping $\cl R$ to $\cl
R^{(k)}$, $k$-copies of $\cl R$, where each $T_i =  M_{\phi_i}$ or 
$M^*_{\phi_i}$ for $\phi_i$ in $H^\infty(\Omega)$ satisfying
$\hat{\phi}_i(\alpha) = 0$.
Suppose that $T$ has closed range. 
Let $\psi$ be a function in $H^\infty(\Omega)$ for which
$\widehat\psi(\alpha)=0$,  and consider the operator 
$T' = \begin{bmatrix}T\\M_\psi \end{bmatrix}$  from $\cl R$ to $\cl
R^{(k+1)}$. Then $T'$ has closed range  and is one-to-one since 
$M_\psi$ is one-to-one. Therefore, the operator 
$T^{'*}T' = \sum_i  T^*_iT_i + M^*_\psi M_\psi$ is invertible.
Hence we obtain a contradiction:
\begin{align*}
 0 \ne \Phi(T'{}^*T')(\alpha)
 &= \sum_i \Phi(T^*_iT_i)(\alpha) + 
     \Phi (M^*_\psi M_\psi)(\alpha)\\
 &= \sum_i |\widehat\phi_i(\alpha)|^2 + |\widehat\psi(\alpha)|^2
  = 0.
\end{align*}

Now suppose $\alpha$ is a point in $M_{H^\infty(\Omega)}$ 
satisfying the statement of the theorem, that is, all operators
of the form 
$T = \begin{sbmatrix}T_1\\ \vdots \\ T_n\end{sbmatrix}$
mapping $\cl  R$ to $\cl R^{(n)}$ with $T_i = M_{\phi_i}$ or 
$M^*_{\phi_i}$ for $\phi_i$ in $H^\infty(\Omega)$ satisfying
$\hat{\phi}_i(\alpha) = 0$ have non-closed range. 
Hence we may construct a net of unit vectors $\{ x_\lambda \}$
indexed by finite subsets of $H^\infty(\Omega)$ as follows.
Let $T_\lambda$ be the column operator with $2|\lambda|$ entries
$T_{\phi - \hat{\phi}(\alpha)}$ and $T_{\phi - \hat{\phi}(\alpha)}^*$
for $\phi \in \lambda$.
Since this is not bounded below, select a unit vector $x_\lambda$ so
that $\| T_\lambda x_\lambda \| < |\lambda|^{-1}$.
By Lemma~\ref{no_eig}, there are no reducing eigenvalues.
So it is easy to verify that the net $x_\lambda$ tends weakly to $0$.

Now following the proof of Theorem~\ref{character}, we can show that
the projections $P_\lambda$ onto the span of $x_\lambda$
asymptotically reduce each $M_\phi$, and hence all of 
$\mf T_\infty( \cl R)$.  It therefore follows as in that proof that
$\Phi(A) = \lim_\lambda (A x_\lambda, x_\lambda)$ is a character
extending evaluation at $\alpha$.
\end{proof}

In \cite{S} Sundberg, using the additional structure present in his 
context, is  able to replace the column operator by a single operator.
In his case, there is a Hilbert super-module $\cl N$ over $A(\Omega)$
containing the given  module $\cl R$ for which all operators $N_\phi$
defined by module multipliers $\phi$ in $A(\Omega)$ (or, equivalently,
in $H^\infty(\Omega)$ if $\Omega$ is weak  pointed) are normal. If
$\mf A_{\cl N}$, respectively $\mf A^\infty_{\cl N}$, is  the
C*-algebra generated by $\{N_\phi : \phi\in A(\Omega)\}$, or 
$\{N_\phi : \phi\in H^\infty(\Omega)\}$, then $\mf A_{\cl N} \simeq 
C(Y_{\cl N})$ and $\mf A^\infty_{\cl N} \simeq C(Y^\infty_{\cl N})$. 
Consideration of the maps from $A(\Omega)$ onto $\mf A_{\cl 
N}$ and $H^\infty(\Omega)$ into $\mf A^\infty_{\cl N}$ allows one to 
identify 
\[
 \partial A(\Omega) \subseteq Y_{\cl N} \subseteq M_{A(\Omega)} \qand 
 \partial H^\infty(\Omega) \subseteq Y^\infty_{\cl N} \subseteq
 M_{H^\infty(\Omega)} .
\]
Let 
$P$ be the projection of $\cl N$ onto $\cl R$. Consideration of the 
matrix  representation for module multipliers in $A(\Omega)$ and
$H^\infty(\Omega)$  relative to the decomposition $\cl N = \cl R
\oplus \cl R^\perp$ enables one to  conclude that $\mf T(\cl R)$ and
$\mf T^\infty(\cl R)$ coincide with the  C*-algebras generated by the
collections $\{T_\phi : \phi\in  C(Y_{\cl N})\}$ and 
$\{T_\phi : \phi\in L^\infty(\mu)\}$, where 
$T_\phi$ is defined to be the ``Toeplitz operator'' $PN_\phi P$  
and $\mu$ is a scalar spectral measure for $\cl N$. Two obvious 
questions which  now present themselves concern the relationship, if
any, between $X_{\cl R}$ and $Y_{\cl R}$  and between $X^\infty_{\cl
R}$ and $Y^\infty_{\cl R}$.

\begin{eg}\label{ex5.4}
Let $\mu$ be the positive measure on the closure of $\bb D$ defined as 
the sum  of Lebesgue measure $d\theta$ on $\partial\bb D$ plus the
atomic measure of  mass one supported at the origin of $\bb D$. If
$H^2(\mu)$ is the closure of  $A(\bb D)$ in $L^2(\mu)$, then it is a
rank one quasi-free Hilbert module over 
$A(\bb D)$. One can see that $Y_{H^2(\mu)}$ is 
$\{0\}\cup \partial\bb D$, while 
$Y^\infty_{H^2(\mu)}$ is $\{0\}\cup M_{L^\infty(d\theta)}$. Moreover, 
one can  show that there are isomorphisms, one of $\mf T(H^2(\bb D))$
with $\mf  T(H^2(\mu))$ and one of $\mf T^\infty(H^2(\bb D))$ with
$\mf  T^\infty(H^2(\mu))$, both preserving the respective inclusion
maps of $A(\bb  D)$ and $H^\infty(\bb D)$ into the C*-algebras.
Therefore, we have 
$X_{H^2(\mu)} = \partial\bb D$ and $X^\infty_{H^2(\mu)} = \partial 
H^\infty(\Omega)$. Since both $X_{H^2(\mu)}\ne Y_{H^2(\mu)}$ and 
$X^\infty_{H^2(\mu)}\ne Y^\infty_{H^2(\mu)}$, we see that neither of 
the pairs 
$X_{\cl R}$ and $Y_{\cl R}$ nor $X^\infty_{\cl R}$ and 
$Y^\infty_{\cl R}$ are equal, in general. However, we can ask whether
inclusion holds, that is, does $X_{\cl R} \subseteq Y_{\cl R}$ or
$X^\infty_{\cl R} \subseteq Y^\infty_{\cl  R}$?
\end{eg}


\end{document}